\documentclass{article}
\usepackage{amsmath}
\usepackage{amssymb}
\usepackage{cite}
\usepackage{tabularx}
\usepackage{graphicx}
\usepackage{graphics}
\usepackage{caption}
\usepackage{subcaption}
\title{{\large\textbf{Convergence and stability of truncated Euler-Maruyama algorithm for stochastic proportional delay Mckean-Vlasov models with jump process}}}

\author{Amr Abosenna$^{1,2}$, Zhuoqi Liu$^{1,*}$\\
{\small\textit{$^1$ Department of Mathematics, Shanghai Normal University, Shanghai 200234, China}}\\
{\small\textit{$^2$ Basic Science Department, Faculty of Engineering at Shoubra, Benha University, 11629, Egypt}}\\
{\small\textit{$^*$ Corresponding author: zhuoqiliu2021@163.com}}}
\date{}

\begin{document}
	\maketitle
	\graphicspath{{Pictures/}}
	\begin{abstract} 
		{\small Stochastic Mckean-Vlasov models have a substantial importance in different fields such as finance, biology and control. This paper puts the light on stochastic proportional delay Mckean-Vlasov model with L\'evy jump where the non-jump coefficients are granted the permission to grow beyond linearity. The truncated Euler-Maruyama algorithm is then applied to our addressed model where the convergence rate and almost sure exponential stability of the aforementioned numerical algorithm are being investigated. Finally, numerical examples are presented to foster the theoretical analysis done throughout the paper.}
	\end{abstract}
	{\small
		\textbf{Keywords:} 	Mckean-Vlasov stochastic proportional delay differential equations, interacting particle system, propagation of chaos, L\'evy jump, truncated Euler Maruyama algorithm, convergence rate, almost sure exponential stability}
	\section{{\normalsize Introduction}}
	{\small
		Stochastic Mckean-Vlasov models have a substantial importance due to their wide utilization in multiple directions such as engineering, biology, physics, mathematical finance and control theory \cite{guhlke2018stochastic,dreyer2011phase,carmona2018probabilistic,bossy2015clarification,bensoussan2013mean,wu2022stabilization}. These models are distribution dependent stochastic differential equations introduced by \cite{mckean1966class} to extend the traditional stochastic differential models by considering the collective behaviour of multiple interacting particles. Furthermore, these types of stochastic models put lights on the statistical distribution of agents rather than individual interactions and include uncertainty in noisy environments, resulting in good control strategies and handling large-scale optimization problems \cite{benachour1998nonlinear,bossy1997stochastic,crisan2010approximate}\\
		
		For achieving realistic simulations for stochastic systems, it is more convenient to take the history data into consideration where the system benefits from the current and past information. These category of stochastic models are called stochastic delay differential equations \cite{babasola2023stochastic,arriojas2007delayed,hobson1998complete} and they have manifested a good performance in stochastic modelling. Stochastic proportional delay differential equations \cite{vivek2022analysis,ren2023stability,ockendon1971dynamics} which are also called stochastic pantograph differential equations are subcategory of stochastic delay differential equations (SDDEs) with vast storage limit and unbounded memory arising in quite different fields of pure and applied mathematics such as number theory, dynamical systems, probability, quantum mechanics and electrodynamics \cite{meng2011pathwise, appleby2016sufficient}. These models can also be regarded as a generalization of pantograph equations. The time delay of these models is called the pantograph delay. The problem arising from a numerical method in dealing with a variable delay is that for the current time step, the delay argument may not hit a previous time step. Thus, there is no previously calculated approximate value of the solution available. Therefore, this problem is tackled by using the approximate value at the nearest grid point on the left of the delay argument. This amounts to interpolation of the undetermined approximate values of the solution at non-mesh points by piecewise constant polynomials whose values are taken at the left end point of the interval containing the delay argument. Therefore, this can be an advantage for using a proportional delay in some cases.\\
		
		Another important issue to be clarified is that Brownian motion is not the best approach for modelling stochastic systems having sudden and considerable changes. Therefore, jump models \cite{merton1976option,kou2002jump,svishchuk2000stochastic,maghsoodi1996mean,higham2005numerical,bruti2007approximation} are the best tools for handling these scenarios due to their crucial role in modelling these sudden changes within the system. Therefore, it is often better to use a jump-diffusion stochastic models for modelling stochastic systems\cite{ bruti2007strong,bass2004stochastic,albeverio2010existence,agrawal2020jump,ahmadian2020exponential}. This paper will address the stochastic Mckean-Vlasov proportional delay model with  L\'evy jump. We refer to \cite{tankov2003financial} for more information regarding L\'evy jump.\\

		Furthermore, most of stochastic Mckean-Vlaso models have no exact solutions and the need for numerical algorithms is indispensable. However, the analysis of these numerical algorithms have been studied under global Lipschitz and linear growth conditions which are somehow restrictive ones \cite{haghighi2019split,milstein2004stochastic}. However, when we try to stray way a little bit beyond these assumptions, problems in convergence analysis appear. This issue was addressed in \cite{hutzenthaler2011strong} and a similar divergence characteristics were noticed in the analysis of stochastic Mckean-Vlasov models \cite{dos2022simulation}. Therefore, attention must be paid during the analysis of numerical algorithms under super linear and non-globally Lipschitz settings. Accordingly, one of the trials to handle this issue was the tamed Euler-Maruyama algorithm which was somehow a novel approach and presented in \cite{hutzenthaler2012strong}. However, due to the modifications done upon the coefficients by the taming technique, the accuracy of the results are not fair enough and this was mentioned in \cite{tretyakov2013fundamental}.\\ 
		
		Another novel approach was implemented by Mao in \cite{mao2015truncated,mao2016convergence} where the truncated Euler-Maruyama technique was addressed and its convergence rate and stability were tackled under somehow relaxed conditions. Guo et al. \cite{guo2018truncated} incorporated delay in the stochastic model and applied the numerical approach discussed in \cite{mao2015truncated} on the stochastic delayed model. This paper will focus on applying the truncated Euler-Maruyama algorithm on addressed model and study the convergence rate and almost sure exponential stability of the numerical algorithm.\\
		
		This work is arranged as follows. Notations, assumptions and model setup are shown in Section 2. Section 3 explains the numerical algorithm. The convergence rate analysis of the numerical algorithm is given in Section 4. Section 5 discusses the almost sure exponential stability of the numerical scheme. Finally, numerical examples are given in Section 6 to elaborate the theoretical analysis.}  
	
	\section{{\normalsize Assumptions \& model setup}}
	{\small
		Let $\big(\Omega,\mathcal{F},(\mathcal{F}_t)_{t\geq0},\mathbb{P}\big)$ be a complete probability space with increasing and right continuous filtration $\left\lbrace \mathcal{F}_{t}\right\rbrace_{t\geq0}$. Let $|\cdot|$ refers to the Euclidean vector norm in $\mathbb R^{d}$ and let $\langle\vartheta,\xi\rangle$ be the inner product of $\vartheta$, $\xi$ in $\mathbb R^{d}$ and for $a\in\mathbb{R}$, $[a]$  extracts only the integer part of $\xi$. Let $\mathcal{P}(\mathbb{R}^d)$ be the set of all probability measures on $\mathbb{R}^d$. For $p\geq1$, define 
		\begin{equation*}
			\mathcal{P}_{p}(\mathbb{R}^d)=\left\lbrace \mu\in\mathcal{P}(\mathbb{R}^d): \mu(|\cdot|^p)=\int_{\mathbb{R}^d}|x|^p\mu(dx)<\infty\right\rbrace.
		\end{equation*} 
		Also, for $p\geq1$, and any $\mu,\nu\in\mathcal{P}_{p}(\mathbb{R}^d)$, the $p-$Wasserstein distance is defined by
		\begin{equation*}
			\mathcal{W}_p(\mu,\nu):=\inf_{\pi^*\in\Pi^{*}(\mu,\nu)}\left(\int_{\mathbb{R}^d\times\mathbb{R}^d}|x-y|^p\pi^*(dx,dy)\right)^{1/p},
		\end{equation*}
		where $\Pi^{*}(\mu,\nu)$ is the family of all couplings for $\mu$ and $\nu$, namely, $\pi^*\in\Pi^{*}(\mu,\nu)$ if and only if $\pi^*(\cdot,\mathbb{R}^d)=\mu(\cdot)$ and $\pi^*(\mathbb{R}^d,\cdot)=\nu(\cdot)$. Clearly, $\mathcal{P}_{p}(\mathbb{R}^d)$ is a Polish space under the $p-$Wasserstein distance $\mathcal{W}_p(\cdot,\cdot)$. Especially, for any $\mu\in\mathcal{P}_{2}(\mathbb{R}^d)$, $\mathcal{W}_2(\mu,\delta_0)=\mu^{\frac{1}{2}}(|\cdot|^2)$, where $\delta_\xi$ is the Dirac measure concentrated at point $\xi\in\mathbb{R}^d$. Let $W(t)$ be $m$-dimensional Brownian motion and $\{\psi(\cdot)\}$ be an $\mathcal{F}_{t}$-adapted L\'evy process with L\'evy measure $\pi(\cdot)$. Also, let $N(\cdot,\cdot)$ be the corresponding  $\mathcal{F}_{t}$-adapted Poisson random measure defined on $R_{+}\times Z$
		\begin{equation*}
			N(t,U):=\sum_{0<s\leq t}I_{U}(\psi(s)-\psi(s^{-})),
		\end{equation*}	
		where $t\geq 0$ and $U$ is a Borel subset of $ Z= R^{d}-\{0\}$. The compensator $\widetilde{N}$ of $N$ is given by 
		\begin{equation*}
			\widetilde{N}(dt,dz)=N(dt,dz)-\pi(dz)dt,
		\end{equation*}		
		where we assume that $W(\cdot)$ and $N(\cdot,\cdot)$ are independent and that $\pi(\cdot)$ is a  L\'evy  measure satisfying 
		\begin{equation*}
			\int_{Z}\mid z\mid^2\pi(dz)<\infty\quad\text{and}\quad 	\pi(Z)<\infty.
		\end{equation*}	 
		Also, let $C$ denotes a general positive constant whose value changes from one place to another.\\\\
		Consider the following stochastic proportional delay Mckean-Vlasov model with L\'evy jump  of the form
		\begin{equation}
			\label{modelmv}
			\begin{aligned}
				dx_{t}&=f(x_{t},x_{\eta t},\upsilon_{t},\upsilon_{\eta t})dt
				+g(x_{t},x_{\eta t},\upsilon_{t},\upsilon_{\eta t})dW_{t}\\
				&\quad+\int_{Z}h(x_{t},x_{\eta t},\upsilon_{t},\upsilon_{\eta t},z)\widetilde{N}(dt,dz),
			\end{aligned}
		\end{equation}
		defined on $ 0\leq t\leq T $ with $0<\eta<1$ and initial data $x_0\in\mathbb{R}^d$ which is $\mathcal{F}_{0}$-measurable random variable with $\mathbb{E}|x_0|^{\tilde{m}}<\infty$ for any $\tilde{m}>0$. Here $\upsilon_{t}$ indicates the law of $x_t$. Also, $f:\mathbb{R}^d \times \mathbb{R}^d \times \mathcal{P}_2(\mathbb{R}^d) \times \mathcal{P}_2(\mathbb{R}^d)\rightarrow \mathbb R^{d}$,  $g:\mathbb{R}^d \times \mathbb{R}^d \times \mathcal{P}_2(\mathbb{R}^d) \times \mathcal{P}_2(\mathbb{R}^d)\rightarrow \mathbb R^{d\times m}$ and $h:\mathbb{R}^d \times \mathbb{R}^d \times \mathcal{P}_2(\mathbb{R}^d) \times \mathcal{P}_2(\mathbb{R}^d) \times Z\rightarrow \mathbb R^{d}, m,d \in \mathbb{N}^{+}$ are Borel measurable functions.\\\\
		\textbf{Assumption 2.1.} \textit{There exist a constant $\varepsilon\geq0$ such that 
			\begin{equation}
				\label{assm1}
				\begin{aligned}
					|f(\vartheta_1,\xi_1&,\nu_1,\nu_1^*)-f(\vartheta_2,\xi_2,\nu_2,\nu_2^*)|\vee
					\|g(\vartheta_1,\xi_1,\nu_1,\nu_1^*)-g(\vartheta_2,\xi_2,\nu_2,\nu_2^*)\|\\
					&\leq C((1+|\vartheta_1|^\varepsilon+|\xi_1|^\varepsilon+|\vartheta_2|^\varepsilon+|\xi_2|^\varepsilon)(|\vartheta_1-\vartheta_2|+|\xi_1-\xi_2|)\\
					&\quad+\mathcal{W}_2(\nu_1,\nu_2)+\mathcal{W}_2(\nu^*_1,\nu^*_2)),
				\end{aligned}
			\end{equation}
			$\text{and}$
			\begin{equation}
				\begin{aligned}
					\int_{Z}|h(\vartheta_1,\xi_1&,\nu_1,\nu_1^*,z)-h(\vartheta_2,\xi_2,\nu_2,\nu_2^*,z)|\pi(dz)\\
					&\leq C(|\vartheta_1-\vartheta_2|+|\xi_1-\xi_2|+\mathcal{W}_2(\nu_1,\nu_2)+\mathcal{W}_2(\nu^*_1,\nu^*_2)),
				\end{aligned}
			\end{equation}
			for all $\vartheta_1,\vartheta_2,\xi_1,\xi_2\in\mathbb{R}^{d}$ and $\nu_1,\nu_2,\nu_1^*,\nu_2^*\in\mathcal{P}_2(\mathbb{R}^d)$. Furthermore
			\begin{equation}
				\mathbb{E}\int_0^T\int_Z\mid h(\vartheta,\xi,\nu,\nu^*,z)\mid^2\pi(dz)dt<\infty,
			\end{equation}
			and
			\begin{equation}
				\int_Z\mid h(\vartheta,\xi,\nu,\nu^*,z)\mid^2\pi(dz)<\infty,
			\end{equation}
			for all $\vartheta,\xi\in\mathbb{R}^{d}$, $\nu,\nu^*\in\mathcal{P}_2(\mathbb{R}^d)$ and $z\in Z$.}\\\\
		\textbf{Assumption 2.2.}\itshape{ There exists a constant $q_0>2$ such that
			\begin{equation}
				\label{assm2}
				\begin{aligned}
					(\vartheta_1&-\vartheta_2)^T(f(\vartheta_1,\xi_1,\nu_1,\nu_1^*)-f(\vartheta_2,\xi_2,\nu_2,\nu_2^*))\\
					&\quad+\frac{q_0-1}{2}
					\|g(\vartheta_1,\xi_1,\nu_1,\nu_1^*)-g(\vartheta_2,\xi_2,\nu_2,\nu_2^*)\|^2\\
					&\leq C(|\vartheta_1-\vartheta_2|^2
					+|\xi_1-\xi_2|^2+\mathcal{W}_2^2(\nu_1,\nu_2)+\mathcal{W}_2^2(\nu^*_1,\nu^*_2)),
				\end{aligned}
			\end{equation}
			for all $\vartheta_1,\vartheta_2,\xi_1,\xi_2\in\mathbb{R}^{d}$ and $\nu_1,\nu_2,\nu_1^*,\nu_2^*\in\mathcal{P}_2(\mathbb{R}^d)$.} \\\\
		\textbf{Assumption 2.3.} \itshape{ There exists a constant $q^*>q_0$ such that
			\begin{equation}
				\label{assm3}
				\begin{aligned}
					(\vartheta)^Tf(\vartheta,\xi,\nu,\nu^*)&+\frac{q^*-1}{2}\|g(\vartheta,\xi,\nu,\nu^*)\|^2\\
					&\leq C(1+|\vartheta|^2
					+|\xi|^2+\mathcal{W}_2^2(\nu,\delta_0)+\mathcal{W}_2^2(\nu^*,\delta_0))
				\end{aligned}
			\end{equation}
			for all $\vartheta,\xi\in\mathbb{R}^{d}$ and $\nu,\nu^*\in\mathcal{P}_2(\mathbb{R}^d)$.}\\\\  
		By exploiting Assumption 2.1, we have 
		\begin{equation}
			\label{ass4}
			\begin{aligned}
				|f(\vartheta,\xi&,\nu,\nu^*)|\vee\|g(\vartheta,\xi,\nu,\nu^*)\|\\
				&\leq C(1+|\vartheta|^{\varepsilon+1}+|\xi|^{\varepsilon+1}|
				+\mathcal{W}_2(\nu,\delta_0)+\mathcal{W}_2(\nu^*,\delta_0)),
			\end{aligned}
		\end{equation}
		and
		\begin{equation}
			\label{assm41}
			\begin{aligned}
				\int_{Z}|h(\vartheta_1,\xi_1,\nu_1,\nu_1^*,z)|\pi(dz)\leq C(1+|\vartheta|+|\xi|+\mathcal{W}_2(\nu,\delta_0)+\mathcal{W}_2(\nu^*,\delta_0)),
			\end{aligned}
		\end{equation}
		for all $\vartheta,\xi\in\mathbb{R}^{d}$ and $\nu,\nu^*\in\mathcal{P}_2(\mathbb{R}^d)$.}\\\\ 
	\textbf{{\small Lemma 2.1.}} \textit{Suppose that Assumptions 2.1-2.3 hold, then Eq.(\ref{modelmv}) has a unique solution and furthermore for any $q\in [2,q^*)$, the solution possesses the property 
		\begin{equation}
			\label{lemma2.1}
			\sup_{0\leq t\leq T}\mathbb{E}|x_t|^{q}\leq C,\quad \forall T>0
	\end{equation}}
	\textbf{{\small Proof.}} The concept of proving this Lemma is similar to the one discussed in \cite{biswas2020well}.\\\\
	The interacting particle system is the common approach to approximate stochastic Mckean-Vlason model, where focal point lies in utilizing the empirical measure to approximate $\upsilon_t$ which is the law of $x_t$. For any $i\in\{1,\cdots,N\}$ where $N\in\mathbb{N}$, let $(x_0^i,W_t^i,\widetilde{N}^i(dt,dz))$ be independent copies of $(x_0,W_t,\widetilde{N}(dt,dz))$. Then, the non-interacting particle system related to Eq.(\ref{modelmv}) is given by
	\begin{equation}
		\begin{aligned}
			dx_{t}^{i}&=f(x_{t}^{i},x_{\eta t}^{i},\upsilon_{t}^{x},\upsilon_{\eta t}^{x})dt
			+g(x_{t}^{i},x_{\eta t}^{i},\upsilon_{t}^{x},\upsilon_{\eta t}^{x})dW_{t}^{i}\\
			&\quad+\int_{Z}h(x_{t}^{i},x_{\eta t}^{i},\upsilon_{t}^{x},\upsilon_{\eta t}^{x},z)\widetilde{N}^i(dt,dz),
		\end{aligned}
	\end{equation}
	with the initial value $x_0^i$, where $\upsilon_{t}^{x}$ is the law of $x_i$ at $t$. To approximate the law, the corresponding interacting particles system is defined by
	\begin{equation}
		\label{modelmv2}
		\begin{aligned}
			dx_{t}^{i,N}&=f(x_{t}^{i,N},x_{\eta t}^{i,N},\upsilon_{t}^{x,N},\upsilon_{\eta t}^{x,N})dt
			+g(x_{t}^{i,N},x_{\eta t}^{i,N},\upsilon_{t}^{x,N},\upsilon_{\eta t}^{x,N})dW_{t}^{i}\\
			&\quad+\int_{Z}h(x_{t}^{i,N},x_{\eta t}^{i,N},\upsilon_{t}^{x,N},\upsilon_{[\eta t]}^{x,N},z)\widetilde{N}^i(dt,dz),
		\end{aligned}
	\end{equation}
	with the initial value $x_0^{i,N}$, where the empirical measure is defined by
	\begin{equation}
		\upsilon_{t}^{x,N}(\cdot):=\frac{1}{N}\sum_{i=1}^{N} \delta_{x_{t}^{i,N}}(\cdot)\quad\text{and}\quad
		\upsilon_{\eta t}^{x,N}(\cdot):=\frac{1}{N}\sum_{i=1}^{N}\delta_{x_{\eta t}^{i,N}}(\cdot).
	\end{equation}
	It is also convenient to write Eq.(\ref{modelmv2}) in stochastic integral form as follows
	\begin{equation}
		\label{integraleq.}
		\begin{aligned}
			x_{t}^{i,N}&=x_{0}^{i,N}+\int_{0}^{t}f(x_{s}^{i,N},x_{\eta s}^{i,N},\upsilon_{s}^{x,N},\upsilon_{\eta s}^{x,N})ds
			+\int_{0}^{t}g(x_{s}^{i,N},x_{\eta s}^{i,N},\upsilon_{s}^{x,N},\upsilon_{\eta s}^{x,N})dW_{s}^{i}\\
			&\quad+\int_{0}^{t}\int_{Z}h(x_{s}^{i,N},x_{\eta s}^{i,N},\upsilon_{s}^{x,N},\upsilon_{\eta s}^{x,N},z)\widetilde{N}^i(ds,dz).
		\end{aligned}
	\end{equation}  
	The following proposition claims that the trajectories of $x_t^{i,N}$ and $x_t^{i}$ come very close to each other as $N\to\infty$ almost surely. This implies the propagation of chaos for the interacting particle system.\\\\
	\textbf{{\small Proposition 2.1.}} \textit{{\small(Propagation of Chaos \cite{guo2024convergence}) Let Assumptions 2.1-2.3 hold.  If for some $p\in [2, q_0/2)$ the initial law has finite $p$-moment, then it holds that
			\begin{equation}
				\label{poc}
				\sup _{i \in\{1, \ldots, N\}} \sup _{t \in[0, T]} \mathbb{E}|x_{t}^{i}-x_{t}^{i,N}|^{p} \leq C
				\left\{\begin{array}{ll}
					N^{-1 / 2}, & \text { if } p>d/2, \\
					N^{-1 / 2} \log(1+N), & \text { if } p=d/2,\\
					N^{-p / d}, & \text { if }p\in[2,d/2),
				\end{array}\right.
			\end{equation}
			where the constant $C> 0$ is independent of $N$.}}

\section{{\normalsize Numerical algorithm}}
{\small
	The drift and diffusion parts may grow super linearly and this can be handled via truncation approach which is implemented by picking up a strictly increasing continuous function $\lambda$ such that
	\begin{equation}
		\sup_{|\vartheta|\vee|\xi|\leq r}\left(|f(\vartheta,\xi,\nu,\nu^*)|\vee\left\|g(\vartheta,\xi,\nu,\nu^*)\right\|\right) \leq \lambda(r),
	\end{equation}	
	where its inverse $\lambda^{-1}$ is strictly increasing continuous function from $[\lambda(0),+\infty)$ to $\mathbb{R}_{+}$. Firstly, the interval $[0,T]$ is divided into $N_{T}\in\mathbb{N}$ subintervals with step size $\varDelta t=T/N_{T}$ where $t_ {l} =l\varDelta t$ for $t_l\in[0,T]$, $l=0,1,\cdots$. Secondly, $\varDelta t^{*}\in (0,1]$ is taken and a strictly decreasing function $\kappa: (0,\varDelta t^{*}]\longrightarrow(0,+\infty)$ is picked up s.t.
	\begin{equation}
		\label{H2condition}
		\kappa\left(\varDelta t^{*}\right) \geq \lambda(1), \quad \lim _{\varDelta t \rightarrow 0} \kappa(\varDelta t)=\infty,\quad(\varDelta t^{1/q^*}\vee\varDelta t^{1/4})\kappa(\varDelta t) <\infty,\quad\forall\varDelta t \in(0,\varDelta t^{*}].
	\end{equation}
	For a given step size $\varDelta t\in (0,\varDelta t^{*}]$, we define the truncated functions by
	\begin{equation}
		\label{trunc1}
		f_{\Delta}(\vartheta,\xi,\nu,\nu^*)=f\biggl((|\vartheta| \wedge \lambda^{-1}(\kappa(\varDelta t)))\frac{\vartheta}{|\vartheta|},\quad
		(|\xi| \wedge \lambda^{-1}(\kappa(\varDelta t)))\frac{\xi}{|\xi|},\quad\nu,\quad \nu^*\biggl),
	\end{equation}
	and
	\begin{equation}
		\label{trunc10}
		g_{\Delta}(\vartheta,\xi,\nu,\nu^*)=g\biggl((|\vartheta| \wedge \lambda^{-1}(\kappa(\varDelta t)))\frac{\vartheta}{|\vartheta|},\quad
		(|\xi| \wedge \lambda^{-1}(\kappa(\varDelta t)))\frac{\xi}{|\xi|},\quad\nu,\quad \nu^*\biggl),
	\end{equation}
	for all $\vartheta,\xi\in\mathbb{R}^{d}$ and $\nu,\nu^*\in\mathcal{P}_2(\mathbb{R}^d)$. By utilizing Assumption 2.3 and following the same approach done in \cite{mao2015truncated}, the truncated coefficients satisfy the following properties 
	\begin{equation}
		\label{Hcondition}
		|f_{\Delta}(\vartheta,\xi,\nu,\nu^*)|\vee\|	g_{\Delta}(\vartheta,\xi,\nu,\nu^*)\|\leq\lambda(\lambda^{-1}(\kappa(\Delta)))=\kappa(\Delta),
	\end{equation}
	and 
	\begin{equation}
		\label{ww}
		\begin{aligned}
			(\vartheta)^Tf_\Delta(\vartheta,\xi,\nu,\nu^*)&+\frac{q^*-1}{2}\|g_\Delta(\vartheta,\xi,\nu,\nu^*)\|^2\\
			&\leq C(1+|\vartheta|^2
			+|\xi|^2+\mathcal{W}_2^2(\nu,\delta_0)+\mathcal{W}_2^2(\nu^*,\delta_0)),
		\end{aligned}
	\end{equation}
	for all $\vartheta,\xi\in\mathbb{R}^{d}$ and $\nu,\nu^*\in\mathcal{P}_2(\mathbb{R}^d)$. Then the numerical approximation of $x_t^{i,N} $ at time $t_l$ is denoted by $\varpi_{l}^{i,N}$. With the defined truncated coefficients given by Eq.(\ref{trunc1}) and Eq. (\ref{trunc10}), the numerical solutions $\varpi_{l}^{i,N}$ are then generated by the basic Euler-Maruyama scheme as follows
	\begin{equation}
		\label{TEM}
		\begin{aligned}	
			\varpi_{l+1}^{i,N}&=\varpi_{l}^{i,N}+ f_{\Delta}(\varpi_{l}^{i,N},\varpi_{[\eta l]}^{i,N},\upsilon_{l}^{\varpi,N},\upsilon_{[\eta l]}^{\varpi,N})\varDelta t
			+g_{\Delta}(\varpi_{l}^{i,N},\varpi_{[\eta l]}^{i,N},\upsilon_{l}^{\varpi,N},\upsilon_{[\eta l]}^{\varpi,N})\Delta W_{l}^{i}\\
			&\quad+\int_{Z}h(\varpi_{l}^{i,N},\varpi_{[\eta l]}^{i,N},\upsilon_{l}^{\varpi,N},\upsilon_{[\eta l]}^{\varpi,N},z)\widetilde{N}^i(\varDelta t,dz),
		\end{aligned}
	\end{equation}
	with initial value $\varpi_{0}^{i,N}=x^{i,N}_0$ where $\Delta W_{l}^{i}=W_{t_{l+1}}^{i}-W_{t_{l}}^{i}$ and
	\begin{equation}
		\upsilon_{l}^{\varpi,N}(\cdot):=\frac{1}{N}\sum_{i=1}^N \delta_{\varpi_{l}^{i,N}}(\cdot)\quad\text{and}\quad
		\upsilon_{[\eta l]}^{\varpi,N}(\cdot):=\frac{1}{N}\sum_{i=1}^N \delta_{\varpi_{[\eta l]}^{i,N}}(\cdot).
	\end{equation} 
	For all $t\in[t_l,t_{l+1})$, we define
	\begin{equation}
		\label{con}
		\begin{aligned}	
			\varpi_{t}^{i,N}&= \varpi_{t_l}^{i,N}+(t-t_l)f_{\Delta}(\varpi_{l}^{i,N},\varpi_{[\eta l]}^{i,N},\upsilon_{l}^{\varpi,N},\upsilon_{[\eta l]}^{\varpi,N})\\
			&\quad+g_{\Delta}(\varpi_{l}^{i,N},\varpi_{[\eta l]}^{i,N},\upsilon_{l}^{\varpi,N},\upsilon_{[\eta l]}^{\varpi,N})(W_{t}^{i}-W_{t_l}^{i})\\
			&\quad+\int_{t_l}^{t}\int_{Z}h(\varpi_{l}^{i,N},\varpi_{[\eta l]}^{i,N},\upsilon_{l}^{\varpi,N},\upsilon_{[\eta l]}^{\varpi,N},z)\widetilde{N}^i(ds,dz),
		\end{aligned}
	\end{equation}
	and denote 
	\begin{equation}
		\begin{aligned}
			{\omega^{*}}_t^{i,N}=\sum_{l=0}^{\infty}\varpi^{i,N}_{l}I_{ [t_{l},t_{l+1})}(t),\quad
			{\tilde\omega}_t^{i,N}=\sum_{l=0}^{\infty}\varpi^{i,N}_{[\eta l]}I_{ [t_{l},t_{l+1})}(t),
		\end{aligned}
	\end{equation}
	and
	\begin{equation}
		\begin{aligned}
			{\upsilon^{*}}_t^{\varpi,N}(\cdot):=\frac{1}{N}\sum_{i=1}^N \delta_{{\omega^{*}}_t^{i,N}}(\cdot),\quad
			{\tilde{\upsilon}}_t^{\varpi,N}(\cdot):=\frac{1}{N}\sum_{i=1}^N \delta_{{\tilde\omega}_t^{i,N}}(\cdot).
		\end{aligned}
	\end{equation}
	Then Eq.(\ref{con}) can be coined in an integral structure as follows
	\begin{equation}
		\label{con2}
		\begin{aligned}	
			\varpi_{t}^{i,N}&=\varpi_0^{i,N}+\int_{0}^{t} f_{\Delta}({\omega^{*}}_s^{i,N},{\tilde\omega}_s^{i,N},{\upsilon^{*}}_s^{\varpi,N},{\tilde{\upsilon}}_s^{\varpi,N})ds\\
			&\quad+\int_{0}^{t}g_{\Delta}({\omega^{*}}_s^{i,N},{\tilde\omega}_s^{i,N},{\upsilon^{*}}_s^{\varpi,N},{\tilde{\upsilon}}_s^{\varpi,N})dW_{s}^{i}\\
			&\quad+\int_{0}^{t}\int_{Z}h({\omega^{*}}_s^{i,N},{\tilde\omega}_s^{i,N},{\upsilon^{*}}_s^{\varpi,N},{\tilde{\upsilon}}_s^{\varpi,N},z)\widetilde{N}^i(ds,dz).
		\end{aligned}
	\end{equation}
\section{{\normalsize Convergence rate analysis}}
{\small
	In this section, we will first set up some lemmas that will help us in proving our main convergence rate result.\\\\
	\textbf{{\small Lemma 4.1.}} {\small\textit{Under Assumption 2.1 and for any $\varDelta t\in (0,\varDelta t^*]$, we have for $p^*\geq2$  
			\begin{equation}
				\label{lemma3.3}
				\begin{aligned}
					\mathbb{E}|&\varpi_{t}^{i,N}-{\omega^{*}}_t^{i,N}|^{p^*}\leq C((\kappa(\varDelta t))^{p^*}(\varDelta t)^{p^*/2}+\varDelta t(1+\mathbb{E}|{\omega^{*}}_t^{i,N}|^{p^*}+\mathbb{E}|{\tilde\omega}_t^{i,N}|^{p^*}))
				\end{aligned}
	\end{equation}}\\
	\textbf{{\small Proof.}} {\small Fix any $\varDelta t\in(0,\varDelta t^*]$, $t\geq0$ and $p^*\geq2$. Then there exists a unique $l$ such that $l\varDelta t\leq t\leq(l+1)\varDelta t$. From Eq. (\ref{con2}), we have the following
		\begin{equation}
			\label{T31}
			\begin{aligned}
				\mathbb{E}|\varpi_{t}^{i,N}-{\omega^{*}}_t^{i,N}|^{p^*}&=\mathbb{E}|\varpi_{t}^{i,N}-\varpi_{l\varDelta t}^{i,N}|^{p^*}\\
				&=\mathbb{E}\Biggl|\int_{l\varDelta t}^{t}f_{\Delta}({\omega^{*}}_s^{i,N},{\tilde\omega}_s^{i,N},{\upsilon^{*}}_s^{\varpi,N},{\tilde{\upsilon}}_s^{\varpi,N})ds\\
				&\quad+\int_{l\varDelta t}^{t}g_{\Delta}({\omega^{*}}_s^{i,N},{\tilde\omega}_s^{i,N},{\upsilon^{*}}_s^{\varpi,N},{\tilde{\upsilon}}_s^{\varpi,N})dW_{s}^{i}\\
				&\quad+\int_{l\varDelta t}^{t}\int_{Z}h({\omega^{*}}_t^{i,N},{\tilde\omega}_t^{i,N},{\upsilon^{*}}_t^{\varpi,N},{\tilde{\upsilon}}_t^{\varpi,N},z)\widetilde{N}^i(dt,dz)\Biggl|^{p^*}.
			\end{aligned}	
	\end{equation}
	Then by utilizing (\ref{Hcondition}), Assumption 2.1 and properties of Itô integral \cite{mao2007stochastic}, we get 
	\begin{equation}
		\label{T32}
		\begin{aligned}
			\mathbb{E}&|\varpi_{t}^{i,N}-{\omega^{*}}_t^{i,N}|^{p^*}\\
			&\leq C(\varDelta t)^{p^*-1}\mathbb{E}\int_{l\varDelta t}^{t}|f_{\Delta}({\omega^{*}}_s^{i,N},{\tilde\omega}_s^{i,N},{\upsilon^{*}}_s^{\varpi,N},{\tilde{\upsilon}}_s^{\varpi,N})|^{p^*}ds\\
			&\quad+C(\varDelta t)^{(p^*-2)/2}\mathbb{E}\int_{l\Delta}^{t}|g_{\Delta}({\omega^{*}}_s^{i,N},{\tilde\omega}_s^{i,N},{\upsilon^{*}}_s^{\varpi,N},{\tilde{\upsilon}}_s^{\varpi,N})|^{p^*}d(s)\\
			&\quad+C\mathbb{E}\Biggl(\int_{l\varDelta t}^{t}\int_{Z}|h({\omega^{*}}_t^{i,N},{\tilde\omega}_t^{i,N},{\upsilon^{*}}_t^{\varpi,N},{\tilde{\upsilon}}_t^{\varpi,N},z)|^{2}\pi(dz)ds\Biggl)^{p^*/2}\\
			&\quad+C\mathbb{E}\int_{l\varDelta t}^{t}\int_{Z}|h({\omega^{*}}_t^{i,N},{\tilde\omega}_t^{i,N},{\upsilon^{*}}_t^{\varpi,N},{\tilde{\upsilon}}_t^{\varpi,N},z)|^{p^*}\pi(dz)ds\\
			&\leq C(\kappa(\varDelta t)^{p^*}(\varDelta t)^{p^{*}/2}\\
			&\quad+\varDelta t(1+\mathbb{E}|{\omega^{*}}_t^{i,N}|^{p^*}+\mathbb{E}|{\tilde\omega}_t^{i,N}|^{p^*}+\mathbb{E}\mathcal{W}_2^{p^*}({\upsilon^{*}}_t^{\varpi,N},\delta_0)+	\mathbb{E}\mathcal{W}_2^{p^*}({\tilde{\upsilon}}_s^{\varpi,N},\delta_0))).
		\end{aligned}	
	\end{equation}
	By utilizing Minkowski inequality and the symmetry among the particles because they are all identically distributed \cite{guo2024convergence}, we have
	\begin{equation}
		\label{minkow}
		\begin{aligned}
			\mathbb{E}\mathcal{W}_2^{p^*}({\upsilon^{*}}_t^{\varpi,N},\delta_0)\leq\mathbb{E}|{\omega^{*}}_t^{i,N}|^{p^*}\quad\text{and}\quad
			\mathbb{E}\mathcal{W}_2^{p^*}({\tilde{\upsilon}}_t^{\varpi,N},\delta_0)\leq\mathbb{E}|{\tilde\omega}_t^{i,N}|^{p^*},\quad \forall i.
		\end{aligned}
	\end{equation}
	Therefore,
	\begin{equation}
		\label{T33_VIP}
		\begin{aligned}
			\mathbb{E}|\varpi_{t}^{i,N}-{\omega^{*}}_t^{i,N}|^{p^*}&\leq 
			C((\kappa(\varDelta t))^{p^*}(\varDelta t)^{p^*/2}+\varDelta t(1+\mathbb{E}|{\omega^{*}}_t^{i,N}|^{p^*}+\mathbb{E}|{\tilde\omega}_t^{i,N}|^{p^*}))
		\end{aligned}	
	\end{equation}
	The proof is complete.\\\\
	\textbf{{\small Corollary 4.1.}} {\small\textit{Suppose that Assumption 2.1 hold. Then for any $\varDelta t\in (0,\varDelta t^*]$ and $t\in[0,T]$, we have for $p^*\geq2$
			\begin{equation}
				\label{lemma3.31}
				\begin{aligned}
					\mathbb{E}|&\varpi_{\eta t}^{i,N}-{\tilde\omega}_t^{i,N}|^{p^*}\leq C((\kappa(\varDelta t))^{p^*}(\varDelta t)^{p^*/2}+\varDelta t(1+\mathbb{E}|{\omega^{*}}_t^{i,N}|^{p^*}+\mathbb{E}|{\tilde\omega}_t^{i,N}|^{p^*}))
				\end{aligned}
	\end{equation}}}\\
	
	\textbf{{\small Proof.}} {\small The proof of this Corollary can be attained by following the same approach as in Lemma 4.1.}\\\\
	\textbf{{\small Lemma 4.2.}} {\small\textit{Let Assumptions 2.1 and 2.3 hold. Then for any $q\in [2,q^*)$, we have  
			\begin{equation}
				\label{lemma3.2}
				\sup_{0<\varDelta t\leq\varDelta t^*}\sup_{0\leq t\leq T}\mathbb{E}|\varpi_{t}^{i,N}|^{q}\leq C,\quad \forall T>0
	\end{equation}}}
	\textbf{{\small Proof.}} {\small For any fixed $\varDelta t\in(0,1]$ and $ t\in[0,T]$, we get by Itô formula \cite{mao2015existence} and (\ref{con2})
		\begin{equation}
			\begin{aligned}
				|\varpi_{t}^{i,N}|^{q}&\leq|x_0^{i,N}|^{q}+\int_{0}^{t}q|\varpi_{s}^{i,N}|^{q-2}\Bigr((\varpi_{s}^{i,N})^{T}f_{\Delta}({\omega^{*}}_s^{i,N},{\tilde\omega}_s^{i,N},{\upsilon^{*}}_s^{\varpi,N},{\tilde{\upsilon}}_s^{\varpi,N})\\
				&\quad+\frac{q-1}{2}\|g_{\Delta}({\omega^{*}}_s^{i,N},{\tilde\omega}_s^{i,N},{\upsilon^{*}}_s^{\varpi,N},{\tilde{\upsilon}}_s^{\varpi,N})\|^{2}\Bigr)ds\\
				&\quad+\int_{0}^{t}q|\varpi_{s}^{i,N}|^{q-2}(\varpi_{s}^{i,N})^{T}g_{\Delta}({\omega^{*}}_s^{i,N},{\tilde\omega}_s^{i,N},{\upsilon^{*}}_s^{\varpi,N},{\tilde{\upsilon}}_s^{\varpi,N})dW^i(s)\\
				&\quad+\int_{0}^{t}\int_{Z}q|\varpi_{s}^{i,N}|^{q-2}{(\varpi_{s}^{i,N})}^Th({\omega^{*}}_t^{i,N},{\tilde\omega}_t^{i,N},{\upsilon^{*}}_t^{\varpi,N},{\tilde{\upsilon}}_t^{\varpi,N},z)\widetilde{N}^i(ds,dz)\\
				&\quad+\int_{0}^{t}\int_{Z}\Bigr[|\varpi_{s}^{i,N}+h({\omega^{*}}_t^{i,N},{\tilde\omega}_t^{i,N},{\upsilon^{*}}_t^{\varpi,N},{\tilde{\upsilon}}_t^{\varpi,N},z)|^{q}-|\varpi_{s}^{i,N}|^{q}\\
				&\quad-q|\varpi_{s}^{i,N}|^{q-2}(\varpi_{s}^{i,N})^{T}h({\omega^{*}}_t^{i,N},{\tilde\omega}_t^{i,N},{\upsilon^{*}}_t^{\varpi,N},{\tilde{\upsilon}}_t^{\varpi,N},z)\Bigr]N^i(ds,dz).\\
			\end{aligned}	
		\end{equation}
		Then
		\begin{equation}
			\begin{aligned}
				|\varpi_{t}^{i,N}|^{q}&\leq|x_0^{i,N}|^{q}+\int_{0}^{t}q|\varpi_{t}^{i,N}|^{q-2}\Bigr(({\omega^{*}}_s^{i,N})^{T}f_{\Delta}({\omega^{*}}_s^{i,N},{\tilde\omega}_s^{i,N},{\upsilon^{*}}_s^{\varpi,N},{\tilde{\upsilon}}_s^{\varpi,N})\\
				&\quad+\frac{q-1}{2}\|g_{\Delta}({\omega^{*}}_s^{i,N},{\tilde\omega}_s^{i,N},{\upsilon^{*}}_s^{\varpi,N},{\tilde{\upsilon}}_s^{\varpi,N})\|^{2}\Bigr)ds\\
				&\quad+\int_{0}^{t}q|\varpi_{s}^{i,N}|^{q-2}(\varpi_{t}^{i,N}-{\omega^{*}}_s^{i,N})^{T}f_{\Delta}({\omega^{*}}_s^{i,N},{\tilde\omega}_s^{i,N},{\upsilon^{*}}_s^{\varpi,N},{\tilde{\upsilon}}_s^{\varpi,N})ds\\
				&\quad+\int_{0}^{t}q|\varpi_{s}^{i,N}|^{q-2}(\varpi_{s}^{i,N})^{T}g_{\Delta}({\omega^{*}}_s^{i,N},{\tilde\omega}_s^{i,N},{\upsilon^{*}}_s^{\varpi,N},{\tilde{\upsilon}}_s^{\varpi,N})dW^i(s)\\
				&\quad+\int_{0}^{t}\int_{Z}q|\varpi_{s}^{i,N}|^{q-2}{(\varpi_{s}^{i,N})}^Th({\omega^{*}}_s^{i,N},{\tilde\omega}_s^{i,N},{\upsilon^{*}}_s^{\varpi,N},{\tilde{\upsilon}}_s^{\varpi,N},z)\widetilde{N}^i(ds,dz)\\
				&\quad+\int_{0}^{t}\int_{Z}\Bigr[|\varpi_{s}^{i,N}+h({\omega^{*}}_s^{i,N},{\tilde\omega}_s^{i,N},{\upsilon^{*}}_s^{\varpi,N},{\tilde{\upsilon}}_s^{\varpi,N},z)|^{q}-|\varpi_{s}^{i,N}|^{q}\\
				&\quad-q|\varpi_{s}^{i,N}|^{q-2}(\varpi_{s}^{i,N})^{T}h({\omega^{*}}_s^{i,N},{\tilde\omega}_s^{i,N},{\upsilon^{*}}_s^{\varpi,N},{\tilde{\upsilon}}_s^{\varpi,N},z)\Bigr]N^i(ds,dz).\\	
			\end{aligned}	
		\end{equation}	
		By applying Assumption 2.3, getting use of Taylor formula \cite{mao2015existence}, Young inequality and then taking Expectation, we get
		\begin{equation}
			\label{main}
			\begin{aligned}
				\mathbb{E}|\varpi_{t}^{i,N}|^{q}&\leq|x_0^{i,N}|^{q}+C\mathbb{E}\int_{0}^{t}|\varpi_{s}^{i,N}|^{q-2}(1+|{\omega^{*}}_s^{i,N}|^{2}+|{\tilde\omega}_s^{i,N}|^{2})ds\\
				&\quad+C\mathbb{E}\int_{0}^{t}|\varpi_{s}^{i,N}|^{q-2}(\varpi_{t}^{i,N}-{\omega^{*}}_s^{i,N})^{T}f_{\Delta}({\omega^{*}}_s^{i,N},{\tilde\omega}_s^{i,N},{\upsilon^{*}}_s^{\varpi,N},{\tilde{\upsilon}}_s^{\varpi,N})ds\\
				&\quad+C\mathbb{E}\int_{0}^{t}\int_{Z}|\varpi_{t}^{i,N}|^{q-2}|h({\omega^{*}}_s^{i,N},{\tilde\omega}_s^{i,N},{\upsilon^{*}}_s^{\varpi,N},{\tilde{\upsilon}}_s^{\varpi,N},z)|^{2}\pi(dz)ds\\
				&\quad+C\mathbb{E}\int_{0}^{t}\int_{Z}|h({\omega^{*}}_s^{i,N},{\tilde\omega}_s^{i,N},{\upsilon^{*}}_s^{\varpi,N},{\tilde{\upsilon}}_s^{\varpi,N},z)|^{q}\pi(dz)ds\\
				&=|x_0|^{q}+L_1+L_2+L_3,
			\end{aligned}	
		\end{equation}
		where
		\begin{equation}
			\label{term1}
			\begin{aligned}
				L_1=C\mathbb{E}\int_{0}^{t}|\varpi_{s}^{i,N}|^{q-2}(1+|{\omega^{*}}_s^{i,N}|^{2}+|{\tilde\omega}_s^{i,N}|^{2})ds,
			\end{aligned}	
		\end{equation}
		\begin{equation}
			\label{term2}
			\begin{aligned}
				L_2=C\mathbb{E}\int_{0}^{t}|\varpi_{s}^{i,N}|^{q-2}(\varpi_{s}^{i,N}-{\omega^{*}}_s^{i,N})^{T}f_{\Delta}({\omega^{*}}_s^{i,N},{\tilde\omega}_s^{i,N},{\upsilon^{*}}_s^{\varpi,N},{\tilde{\upsilon}}_s^{\varpi,N})ds,
			\end{aligned}	
		\end{equation}
		and
		\begin{equation}
			\label{term4}
			\begin{aligned}
				L_3&=C\mathbb{E}\int_{0}^{t}\int_{Z}\Bigr(|\varpi_{t}^{i,N}|^{q-2}|h({\omega^{*}}_s^{i,N},{\tilde\omega}_s^{i,N},{\upsilon^{*}}_s^{\varpi,N},{\tilde{\upsilon}}_s^{\varpi,N},z)|^{2}\\
				&\quad+|h({\omega^{*}}_s^{i,N},{\tilde\omega}_s^{i,N},{\upsilon^{*}}_s^{\varpi,N},{\tilde{\upsilon}}_s^{\varpi,N},z)|^{q}\Bigr)\pi(dz)ds.\\
			\end{aligned}	
	\end{equation}}
	
	{\small By utilizing Young inequality, we get
		\begin{equation}
			\label{T1}
			\begin{aligned}
				L_1&\leq C\mathbb{E}\int_{0}^{t}|\varpi_{s}^{i,N}|^{q}ds+C\mathbb{E}\int_{0}^{t}(1+|{\omega^{*}}_s^{i,N}|^{q}+|{\tilde\omega}_s^{i,N}|^{q})ds\\
				&\leq C+C\int_{0}^{t}(\mathbb{E}|\varpi_{s}^{i,N}|^{q}+\mathbb{E}|{\omega^{*}}_s^{i,N}|^{q}+\mathbb{E}|{\tilde\omega}_s^{i,N}|^{q})ds\\
				&\leq C+C\int_{0}^{t}\sup_{0\leq r\leq s}\mathbb{E}|\varpi_{r}^{i,N}|^{q}ds
			\end{aligned}
		\end{equation}	
		By applying Young inequality, Lemma 4.1, (\ref{H2condition}) and (\ref{Hcondition}),we get
		\begin{equation}
			\label{T2}
			\begin{aligned}
				L_2&\leq C+C\mathbb{E}\int_{0}^{t}|\varpi_{s}^{i,N}|^{q}ds
				+C(\kappa(\varDelta t))^{q}\int_{0}^{t}((\kappa(\varDelta t))^{{q}}(\varDelta t)^{{q}/2}\\
				&\quad+\varDelta t(1+\mathbb{E}|{\omega^{*}}_s^{i,N}|^{q}+\mathbb{E}|{\tilde\omega}_s^{i,N}|^{q}))ds\\
				&\leq C+C\int_{0}^{t}\sup_{0\leq r\leq s}\mathbb{E}|\varpi_{r}^{i,N}|^{q}ds
			\end{aligned}	
		\end{equation}
		By utilizing Young inequality and Assumption 2.1 then proceeding the same as it was done before, we get 
		\begin{equation}
			\label{T4}
			\begin{aligned}
				L_3&\leq C+C\int_{0}^{t}\sup_{0\leq r\leq s}\mathbb{E}|y_{\Delta}(r)|^{q}ds.
			\end{aligned}
		\end{equation}
		By plugging in (\ref{T1}), (\ref{T2}) and (\ref{T4}) into (\ref{main}), we get
		\begin{equation}
			\label{tyt}
			\mathbb{E}|\varpi_{t}^{i,N}|^{q}\leq C+C\int_{0}^{t}\sup_{0\leq r\leq s}\mathbb{E}|\varpi_{r}^{i,N}|^{q}ds,
		\end{equation}	
		where the right-hand side of (\ref{tyt}) is increasing in $t$. Then we have
		\begin{equation}
			\sup_{0\leq r\leq t}\mathbb{E}|\varpi_{r}^{i,N}|^{q}\leq C+C\int_{0}^{t}\sup_{0\leq r\leq s}\mathbb{E}|\varpi_{r}^{i,N}|^{q}ds.
		\end{equation}
		Applying the Gronwall inequality yields
		\begin{equation}
			\sup_{0\leq r\leq t}\mathbb{E}|\varpi_{r}^{i,N}|^{q}\leq C.
		\end{equation}
		As this holds for any $ 0<\varDelta t\leq1$ while $ C $ is independent of $\varDelta t $, the required assertion (\ref{lemma3.2}) is obtained.}\\\\
	\textbf{{\small Lemma 4.3.}}{\small\textit{ Let Assumptions 2.1 and 2.3 hold. Then for any $\varDelta t\in (0,\varDelta t^*]$ and $ t\in[0,T] $, we have  
			\begin{equation}
				\label{3a}
				\begin{aligned}
					\mathbb{E}|\varpi_{t}^{i,N}-{\omega^{*}}_t^{i,N}|^{q}&\leq 
					C((\kappa(\varDelta t))^{q}(\varDelta t)^{q/2}+\varDelta t),\quad q\geq2.	
				\end{aligned}
			\end{equation}
			\begin{equation}
				\label{3b}
				\begin{aligned}	
					\mathbb{E}|&\varpi_{\eta t}^{i,N}-{\tilde\omega}_t^{i,N}|^{q}\leq C((\kappa(\varDelta t))^{q}(\varDelta t)^{q/2}+\varDelta t),\quad q\geq2.
				\end{aligned}
	\end{equation}}}
	\textbf{{\small Proof.}} {\small By utilizing Lemmas 4.1 and 4.2 and Corollary 4.1, the required assertions (\ref{3a}) and (\ref{3b}) are directly attained and the proof is complete.}\\\\
	\textbf{{\small Lemma 4.4.}} {\small\textit{Let Assumptions 2.1 and 2.3 hold. For any real number $v>|x_0^{i,N}|$ and $\varDelta t\in (0,\varDelta t^*]$, we define the stopping time $ \tau_v^i=\inf\{t\geq 0:|x_t^{i,N}|\geq v\}$ such that  
			\begin{equation}
				\label{stopping1}
				\mathbb{P}(\tau_v^i\leq T)\leq \frac{C}{v^{q}}
	\end{equation}}}
	\textbf{{\small Proof.}} {\small By utilizing (\ref{lemma2.1}), we have
		\begin{equation*}
			\label{sto1}
			\sup_{0\leq s\leq T}\mathbb{E}|x_{s\wedge\tau_v^i}^{i,N}|^{q}\leq C.
		\end{equation*}
		Then by using applying Chebyshev's inequality, we have
		\begin{equation}
			v^{q}\mathbb{P}(\tau_v^i\leq T)\leq C.
		\end{equation} The proof is done.}\\\\
	\textbf{{\small Lemma 4.5.}} {\small\textit{Let Assumptions 2.1 and 2.3 hold. For any real number $u>|\varpi_{0}^{i,N}|$ and $\varDelta t\in (0,\varDelta t^*]$, we define stopping times $ \tau^{*,i}_{u}=\inf\{t\geq 0:|\varpi_{t}^{i,N}|\geq u\}$ such that  
			\begin{equation}
				\label{stopping2}
				\mathbb{P}(\tau^{*,i}_{u}\leq T)\leq \frac{C}{u^{q}}
	\end{equation}}}
	\textbf{{\small Proof.}} {\small By proceeding in the same manner as it was done in Lemma 4.2, it can be shown that 
		\begin{equation*}
			\label{sto2}
			\sup_{0\leq s\leq T}\mathbb{E}|\varpi_{s\wedge\tau^{*,i}_{u}}^{i,N}|^{q}\leq C.
		\end{equation*}
		Then by applying Chebyshev's inequality, we get
		\begin{equation}
			u^{q}\mathbb{P}(\tau^{*,i}_{u}\leq T)\leq C.
		\end{equation}
		The proof is complete.}\\\\
	\textbf{{\small Theorem 4.1.}} {\small\textit{Let Assumptions 2.1-2.3 hold. Suppose that there exists a real number $q\in(2,q^*)$ such that $q>(1+\varepsilon)q_0$. Then, for any $p\in[2,q_0)  $ and $\varDelta t\in(0,\varDelta t^*]$
			\begin{equation}
				\label{main66}
				\begin{aligned}
					\sup _{t\in[0, T]}\mathbb{E}|x_{t}^{i}-\varpi_{t}^{i,N}|^{p}
					\leq C
					\left\{\begin{array}{ll}
						N^{-1 / 2}+D, & \text { if } p>d/2, \\
						N^{-1 / 2} \log(1+N)+D, & \text { if } p=d/2,\\
						N^{-p / d}+D, & \text { if }p\in[2,d/2),
					\end{array}\right.
				\end{aligned}
			\end{equation}
			where
			\begin{equation}
				\label{main02}
				D=C((\lambda^{-1}(\kappa(\varDelta t)))^{p\varepsilon+p-q}+(\kappa(\varDelta t))^{p}\varDelta t^{p/2}+\varDelta t^{(q-p\varepsilon)/q}).
			\end{equation}
			In particular, we define
			\begin{equation}
				\label{qq}
				\lambda(r)=Cr^{1+\varepsilon}, \quad r\geq 0,
			\end{equation}	
			and let
			\begin{equation}
				\label{qqq}
				\kappa(\varDelta t)=\varDelta t^{-\gamma}, \quad\text{for some}\quad \gamma\in (0,1/4\wedge 1/q].
			\end{equation}
			Then, we get 
			\begin{equation}
				\label{main660}
				\begin{aligned}
					\sup _{t\in[0, T]}\mathbb{E}|x_{t}^{i}-\varpi_{t}^{i,N}|^{p}
					\leq C
					\left\{\begin{array}{ll}
						N^{-1 / 2}+D^*, & \text { if } p>d/2, \\
						N^{-1 / 2} \log(1+N)+D^*, & \text { if } p=d/2,\\
						N^{-p / d}+D^*, & \text { if }p\in[2,d/2),
					\end{array}\right.
				\end{aligned}
			\end{equation}
			where
			\begin{equation}
				\label{main20}
				D^*=C\varDelta t^{[\gamma(q-(1+\varepsilon)p)/(1+\varepsilon)]\wedge[p(1-2\gamma)/2]\wedge[(q-\varepsilon p)/q]}
			\end{equation}
			for all $\varDelta t\in(0,\varDelta t^*]$.
	}}\\\\			
	\textbf{{\small Proof.}} {\small Let $e_i(t)=x_{t}^{i,N}-\varpi_{t}^{i,N}$ and $\tau_i=\tau_v^i\wedge\tau^{*,i}_{u}$. Then by applying Itô formula, we have
		\begin{equation}
			\label{error1}
			\begin{aligned}
				|e_i&(t\wedge\tau_i)|^{p}\\
				&\leq \int_{0}^{t\wedge\tau_i}p|e_i(s)|^{p-2}\Bigr(e_i^{T}(s)[f(x_{s}^{i,N},x_{[\eta s]}^{i,N},\upsilon_{s}^{x,N},\upsilon_{[\eta s]}^{x,N})\\
				&\quad-f_{\Delta}({\omega^{*}}_s^{i,N},{\tilde\omega}_s^{i,N},{\upsilon^{*}}_s^{\varpi,N},{\tilde{\upsilon}}_s^{\varpi,N})]\\
				&\quad+\frac{p-1}{2}|g(x_{s}^{i,N},x_{[\eta s]}^{i,N},\upsilon_{s}^{x,N},\upsilon_{[\eta s]}^{x,N})\\
				&\quad-g_{\Delta}({\omega^{*}}_s^{i,N},{\tilde\omega}_s^{i,N},{\upsilon^{*}}_s^{\varpi,N},{\tilde{\upsilon}}_s^{\varpi,N})|^{2}\Bigr)ds\\
				&\quad+\int_{0}^{t\wedge\tau_i}p|e_{i}(s)|^{p-2}e_{i}^{T}(s)[g(x_{s}^{i,N},x_{[\eta s]}^{i,N},\upsilon_{s}^{x,N},\upsilon_{[\eta s]}^{x,N})\\
				&\quad-g_{\Delta}({\omega^{*}}_s^{i,N},{\tilde\omega}_s^{i,N},{\upsilon^{*}}_s^{\varpi,N},{\tilde{\upsilon}}_s^{\varpi,N})]dW^{i}(s)\\
				&\quad+\int_{0}^{t\wedge\tau_i}\int_{Z}p|e_i(s)|^{p-2}e_i^{T}(s)[h(x_{s}^{i,N},x_{[\eta s]}^{i,N},\upsilon_{s}^{x,N},\upsilon_{[\eta s]}^{x,N},z)\\
				&\quad-h({\omega^{*}}_s^{i,N},{\tilde\omega}_s^{i,N},{\upsilon^{*}}_s^{\varpi,N},{\tilde{\upsilon}}_s^{\varpi,N},z)]\widetilde{N}^{i}(ds,dz)\\
				&\quad+\int_{0}^{t\wedge\tau_i}\int_{Z}\Bigr[|e_i(s)+h(x_{s}^{i,N},x_{[\eta s]}^{i,N},\upsilon_{s}^{x,N},\upsilon_{[\eta s]}^{x,N},z)\\
				&\quad-h({\omega^{*}}_s^{i,N},{\tilde\omega}_s^{i,N},{\upsilon^{*}}_s^{\varpi,N},{\tilde{\upsilon}}_s^{\varpi,N},z)|^{p}-|e_i(s)|^{p}\\
				&\quad-p|e_i(s)|^{p-2}e_i^{T}(s)[h(x_{s}^{i,N},x_{[\eta s]}^{i,N},\upsilon_{s}^{x,N},\upsilon_{[\eta s]}^{x,N},z)\\
				&\quad-h({\omega^{*}}_s^{i,N},{\tilde\omega}_s^{i,N},{\upsilon^{*}}_s^{\varpi,N},{\tilde{\upsilon}}_s^{\varpi,N},z)]\Bigr]N^{i}(ds,dz).\\
			\end{aligned}	
		\end{equation}
		Then by getting use of Taylor formula and taking Expectation, we obtain
		\begin{equation}
			\label{error2}
			\begin{aligned}
				\mathbb{E}&|e_{i}(t\wedge\tau_i)|^{p}\\
				&\leq C\mathbb{E}\int_{0}^{t\wedge\tau_i}p|e_{i}(s)|^{p-2}\\
				&\quad\times\Bigr(e_{i}^{T}(s)[f(x_{s}^{i,N},x_{[\eta s]}^{i,N},\upsilon_{s}^{x,N},\upsilon_{[\eta s]}^{x,N})
				-f_{\Delta}({\omega^{*}}_s^{i,N},{\tilde\omega}_s^{i,N},{\upsilon^{*}}_s^{\varpi,N},{\tilde{\upsilon}}_s^{\varpi,N})]\\
				&\quad+\frac{p-1}{2}|g(x_{s}^{i,N},x_{[\eta s]}^{i,N},\upsilon_{s}^{x,N},\upsilon_{[\eta s]}^{x,N})
				-g_{\Delta}({\omega^{*}}_s^{i,N},{\tilde\omega}_s^{i,N},{\upsilon^{*}}_s^{\varpi,N},{\tilde{\upsilon}}_s^{\varpi,N})|^{2}\Bigr)ds\\
				&\quad+C\mathbb{E}\int_{0}^{t\wedge\tau_i}\int_{Z}\Bigr(|e_{i}(s)|^{p-2}\\
				&\quad\times|h(x_{s}^{i,N},x_{[\eta s]}^{i,N},\upsilon_{s}^{x,N},\upsilon_{[\eta s]}^{x,N},z)
				-h({\omega^{*}}_s^{i,N},{\tilde\omega}_s^{i,N},{\upsilon^{*}}_s^{\varpi,N},{\tilde{\upsilon}}_s^{\varpi,N},z)|^{2}\\
				&\quad+|h(x_{s}^{i,N},x_{[\eta s]}^{i,N},\upsilon_{s}^{x,N},\upsilon_{[\eta s]}^{x,N},z)
				-h({\omega^{*}}_s^{i,N},{\tilde\omega}_s^{i,N},{\upsilon^{*}}_s^{\varpi,N},{\tilde{\upsilon}}_s^{\varpi,N},z)|^{p}\Bigr)\pi(dz)ds.\\
			\end{aligned}	
		\end{equation}
		Furthermore it can be deduced that for $m^*\in(p,q_0)$, we have $q>(1+\varepsilon)m^*$ and
		\begin{equation}
			\label{three}
			\begin{aligned}
				\frac{p-1}{2}&|g(x_{s}^{i,N},x_{[\eta s]}^{i,N},\upsilon_{s}^{x,N},\upsilon_{[\eta s]}^{x,N})
				-g_{\Delta}({\omega^{*}}_s^{i,N},{\tilde\omega}_s^{i,N},{\upsilon^{*}}_s^{\varpi,N},{\tilde{\upsilon}}_s^{\varpi,N})|^{2}\\
				&\leq C(|g(x_{s}^{i,N},x_{[\eta s]}^{i,N},\upsilon_{s}^{x,N},\upsilon_{[\eta s]}^{x,N})-g(\varpi_{s}^{i,N},\varpi_{[\eta s]}^{i,N},\upsilon_{s}^{\varpi,N},\upsilon_{[\eta s]}^{\varpi,N})|^{2}\\
				&\quad+|g(\varpi_{s}^{i,N},\varpi_{[\eta s]}^{i,N},\upsilon_{s}^{\varpi,N},\upsilon_{[\eta s]}^{\varpi,N})-g_{\Delta}({\omega^{*}}_s^{i,N},{\tilde\omega}_s^{i,N},{\upsilon^{*}}_s^{\varpi,N},{\tilde{\upsilon}}_s^{\varpi,N})|^{2})\\
				&= C(|g(x_{s}^{i,N},x_{[\eta s]}^{i,N},\upsilon_{s}^{x,N},\upsilon_{[\eta s]}^{x,N})-g(\varpi_{s}^{i,N},\varpi_{[\eta s]}^{i,N},\upsilon_{s}^{\varpi,N},\upsilon_{[\eta s]}^{\varpi,N})|^{2}\\
				&\quad+|g(\varpi_{s}^{i,N},\varpi_{[\eta s]}^{i,N},\upsilon_{s}^{\varpi,N},\upsilon_{[\eta s]}^{\varpi,N})-g_{\Delta}({\omega^{*}}_s^{i,N},{\tilde\omega}_s^{i,N},{\upsilon^{*}}_s^{\varpi,N},{\tilde{\upsilon}}_s^{\varpi,N})|^{2}).	
			\end{aligned}	
		\end{equation}	 
		Plugging (\ref{three}) into (\ref{error2}) yields
		\begin{equation}
			\label{error4}
			\begin{aligned}
				\mathbb{E}|e_{i}(t\wedge\tau_i)|^{p}
				&\leq C\mathbb{E}\int_{0}^{t\wedge\tau_i}p|e_{i}(s)|^{p-2}
				\Bigr(e_{i}^{T}(s)[f(x_{s}^{i,N},x_{[\eta s]}^{i,N},\upsilon_{s}^{x,N},\upsilon_{[\eta s]}^{x,N})\\
				&\quad-f(\varpi_{s}^{i,N},\varpi_{[\eta s]}^{i,N},\upsilon_{s}^{\varpi,N},\upsilon_{[\eta s]}^{\varpi,N})]
				+\frac{m^*-1}{2}|g(x_{s}^{i,N},x_{[\eta s]}^{i,N},\upsilon_{s}^{x,N},\upsilon_{[\eta s]}^{x,N})\\
				&\quad-g(\varpi_{s}^{i,N},\varpi_{[\eta s]}^{i,N},\upsilon_{s}^{\varpi,N},\upsilon_{[\eta s]}^{\varpi,N})|^{2}\Bigr)ds\\
				&\quad+C\mathbb{E}\int_{0}^{t\wedge\tau_i}p|e_{i}(s)|^{p-2}
				\Bigr(e_{i}^{T}(s)[f(\varpi_{s}^{i,N},\varpi_{[\eta s]}^{i,N},\upsilon_{s}^{\varpi,N},\upsilon_{[\eta s]}^{\varpi,N})\\
				&\quad-f_{\Delta}({\omega^{*}}_s^{i,N},{\tilde\omega}_s^{i,N},{\upsilon^{*}}_s^{\varpi,N},{\tilde{\upsilon}}_s^{\varpi,N})]
				+|g(\varpi_{s}^{i,N},\varpi_{[\eta s]}^{i,N},\upsilon_{s}^{\varpi,N},\upsilon_{[\eta s]}^{\varpi,N})\\
				&\quad-g_{\Delta}({\omega^{*}}_s^{i,N},{\tilde\omega}_s^{i,N},{\upsilon^{*}}_s^{\varpi,N},{\tilde{\upsilon}}_s^{\varpi,N})|^{2}\Bigr)ds\\
				&\quad+C\mathbb{E}\int_{0}^{t\wedge\tau_i}\int_{Z}\Bigr(|e_{i}(s)|^{p-2}
				|h(x_{s}^{i,N},x_{[\eta s]}^{i,N},\upsilon_{s}^{x,N},\upsilon_{[\eta s]}^{x,N},z)\\
				&\quad-h({\omega^{*}}_s^{i,N},{\tilde\omega}_s^{i,N},{\upsilon^{*}}_s^{\varpi,N},{\tilde{\upsilon}}_s^{\varpi,N},z)|^{2}
				+|h(x_{s}^{i,N},x_{[\eta s]}^{i,N},\upsilon_{s}^{x,N},\upsilon_{[\eta s]}^{x,N},z)\\
				&\quad-h({\omega^{*}}_s^{i,N},{\tilde\omega}_s^{i,N},{\upsilon^{*}}_s^{\varpi,N},{\tilde{\upsilon}}_s^{\varpi,N},z)|^{p}\Bigr)\pi(dz)ds.\\
				&=J_1+J_2+J_3,
			\end{aligned}	
		\end{equation}
		where
		\begin{equation}
			\label{j1}
			\begin{aligned}
				J_1&=C\mathbb{E}\int_{0}^{t\wedge\tau_i}p|e_{i}(s)|^{p-2}
				\Bigr(e_{i}^{T}(s)[f(x_{s}^{i,N},x_{[\eta s]}^{i,N},\upsilon_{s}^{x,N},\upsilon_{[\eta s]}^{x,N})\\
				&\quad-f(\varpi_{s}^{i,N},\varpi_{[\eta s]}^{i,N},\upsilon_{s}^{\varpi,N},\upsilon_{[\eta s]}^{\varpi,N})]
				+\frac{m^*-1}{2}|g(x_{s}^{i,N},x_{[\eta s]}^{i,N},\upsilon_{s}^{x,N},\upsilon_{[\eta s]}^{x,N})\\
				&\quad-g(\varpi_{s}^{i,N},\varpi_{[\eta s]}^{i,N},\upsilon_{s}^{\varpi,N},\upsilon_{[\eta s]}^{\varpi,N})|^{2}\Bigr)ds,
			\end{aligned}	
		\end{equation}	
		\begin{equation}
			\label{j2}
			\begin{aligned}
				J_2&=C\mathbb{E}\int_{0}^{t\wedge\tau_i}p|e_{i}(s)|^{p-2}
				\Bigr(e_{i}^{T}(s)[f(\varpi_{s}^{i,N},\varpi_{[\eta s]}^{i,N},\upsilon_{s}^{\varpi,N},\upsilon_{[\eta s]}^{\varpi,N})\\
				&\quad-f_{\Delta}({\omega^{*}}_s^{i,N},{\tilde\omega}_s^{i,N},{\upsilon^{*}}_s^{\varpi,N},{\tilde{\upsilon}}_s^{\varpi,N})]
				+|g(\varpi_{s}^{i,N},\varpi_{[\eta s]}^{i,N},\upsilon_{s}^{\varpi,N},\upsilon_{[\eta s]}^{\varpi,N})\\
				&\quad-g_{\Delta}({\omega^{*}}_s^{i,N},{\tilde\omega}_s^{i,N},{\upsilon^{*}}_s^{\varpi,N},{\tilde{\upsilon}}_s^{\varpi,N})|^{2}\Bigr)ds,
			\end{aligned}	
		\end{equation}
		and
		\begin{equation}
			\label{j3}
			\begin{aligned}
				J_3&=C\mathbb{E}\int_{0}^{t\wedge\tau_i}\int_{Z}\Bigr(|e_{i}(s)|^{p-2}
				|h(x_{s}^{i,N},x_{[\eta s]}^{i,N},\upsilon_{s}^{x,N},\upsilon_{[\eta s]}^{x,N},z)\\
				&\quad-h({\omega^{*}}_s^{i,N},{\tilde\omega}_s^{i,N},{\upsilon^{*}}_s^{\varpi,N},{\tilde{\upsilon}}_s^{\varpi,N},z)|^{2}
				+|h(x_{s}^{i,N},x_{[\eta s]}^{i,N},\upsilon_{s}^{x,N},\upsilon_{[\eta s]}^{x,N},z)\\
				&\quad-h({\omega^{*}}_s^{i,N},{\tilde\omega}_s^{i,N},{\upsilon^{*}}_s^{\varpi,N},{\tilde{\upsilon}}_s^{\varpi,N},z)|^{p}\Bigr)\pi(dz)ds.
			\end{aligned}	
		\end{equation}
		By utilizing Assumption 2.2, it can directly concluded that 
		\begin{equation}
			\label{j11}
			J_1\leq C\int_{0}^{t}\mathbb{E}|e_{i}(s\wedge\tau_i)|^{p}ds.
		\end{equation}	
		\begin{equation}
			\label{j22}
			\begin{aligned}
				J_2&\leq C\mathbb{E}\int_{0}^{t\wedge\tau_i}p|e_{i}(s)|^{p-2}\\
				&\quad\times\Bigr(e_{i}^{T}(s)[f(\varpi_{s}^{i,N},\varpi_{[\eta s]}^{i,N},\upsilon_{s}^{\varpi,N},\upsilon_{[\eta s]}^{\varpi,N})
				-f_{\Delta}(\varpi_{s}^{i,N},\varpi_{[\eta s]}^{i,N},\upsilon_{s}^{\varpi,N},\upsilon_{[\eta s]}^{\varpi,N})]\\
				&\quad+|g(\varpi_{s}^{i,N},\varpi_{[\eta s]}^{i,N},\upsilon_{s}^{\varpi,N},\upsilon_{[\eta s]}^{\varpi,N})
				-g_{\Delta}(\varpi_{s}^{i,N},\varpi_{[\eta s]}^{i,N},\upsilon_{s}^{\varpi,N},\upsilon_{[\eta s]}^{\varpi,N})|^{2}\Bigr)ds\\
				&\quad+C\mathbb{E}\int_{0}^{t\wedge\tau_i}p|e_{i}(s)|^{p-2}\\
				&\quad\times\Bigr(e_{i}^{T}(s)[f_{\Delta}(\varpi_{s}^{i,N},\varpi_{[\eta s]}^{i,N},\upsilon_{s}^{\varpi,N},\upsilon_{[\eta s]}^{\varpi,N})-f_{\Delta}({\omega^{*}}_s^{i,N},{\tilde\omega}_s^{i,N},{\upsilon^{*}}_s^{\varpi,N},{\tilde{\upsilon}}_s^{\varpi,N})]\\
				&\quad+|g_{\Delta}(\varpi_{s}^{i,N},\varpi_{[\eta s]}^{i,N},\upsilon_{s}^{\varpi,N},\upsilon_{[\eta s]}^{\varpi,N})-g_{\Delta}({\omega^{*}}_s^{i,N},{\tilde\omega}_s^{i,N},{\upsilon^{*}}_s^{\varpi,N},{\tilde{\upsilon}}_s^{\varpi,N})|^{2}\Bigr)ds\\
				&=J_{21}+J_{22}.
			\end{aligned}	
		\end{equation}
		By utilizing Young inequality, H\"{o}lder inequality, (\ref{trunc1}),  Assumption 2.1 and (\ref{trunc2}), we get
		\begin{equation}
			\label{j221}
			\begin{aligned}
				J_{21}&\leq C\mathbb{E}\int_{0}^{t}|e_{i}(s\wedge\tau_i)|^{p}ds\\
				&\quad+C\mathbb{E}\int_{0}^{T}\Bigr(|f(\varpi_{s}^{i,N},\varpi_{[\eta s]}^{i,N},\upsilon_{s}^{\varpi,N},\upsilon_{[\eta s]}^{\varpi,N})
				-f_{\Delta}(\varpi_{s}^{i,N},\varpi_{[\eta s]}^{i,N},\upsilon_{s}^{\varpi,N},\upsilon_{[\eta s]}^{\varpi,N})|^{p}\\
				&\quad+|g(\varpi_{s}^{i,N},\varpi_{[\eta s]}^{i,N},\upsilon_{s}^{\varpi,N},\upsilon_{[\eta s]}^{\varpi,N})
				-g_{\Delta}(\varpi_{s}^{i,N},\varpi_{[\eta s]}^{i,N},\upsilon_{s}^{\varpi,N},\upsilon_{[\eta s]}^{\varpi,N})|^{p}\Bigr)ds\\
				&\leq C\int_{0}^{t}\mathbb{E}|e_{i}(s\wedge\tau_i)|^{p}ds+C\mathbb{E}\int_{0}^{T}\Bigr[(1+|\varpi_{s}^{i,N}|^{p\varepsilon}+|\varpi_{[\eta s]}^{i,N}|^{p\varepsilon}
				+|\Gamma(\varpi_{s}^{i,N})|^{p\varepsilon}\\
				&\quad+|\Gamma(\varpi_{[\eta s]}^{i,N})|^{p\varepsilon})
				\times(|\varpi_{s}^{i,N}-\Gamma(\varpi_{s}^{i,N})|^{p}
				+|\varpi_{[\eta s]}^{i,N}-\Gamma(\varpi_{[\eta s]}^{i,N})|^{p})\Bigr]ds\\
				&\leq C\int_{0}^{t}\mathbb{E}|e_{i}(s\wedge\tau_i)|^{p}ds+C\int_{0}^{T}\Bigr[\mathbb{E}(1+|\varpi_{s}^{i,N}|^{q}+|\varpi_{[\eta s]}^{i,N}|^{q}\Bigr]^{p\varepsilon/q}\\
				&\quad\times\Bigr[\mathbb{E}|\varpi_{s}^{i,N}-\Gamma(\varpi_{s}^{i,N})|^{{pq/(q-p\varepsilon)}}+\mathbb{E}|\varpi_{[\eta s]}^{i,N}-\Gamma(\varpi_{[\eta s]}^{i,N})|^{pq/(q-p\varepsilon)}\Bigr]^{(q-p\varepsilon)/q}ds.
			\end{aligned}
		\end{equation}
		Then by utilizing Lemma 4.2, we have
		\begin{equation}
			\label{j221_}
			\begin{aligned}
				J_{21}&\leq C\int_{0}^{t}\mathbb{E}|e_{i}(s\wedge\tau_i)|^{p}ds+C\int_{0}^{T}[\mathbb{E}|\varpi_{s}^{i,N}-\Gamma(\varpi_{s}^{i,N})|^{{pq/(q-p\varepsilon)}}]^{(q-p\varepsilon)/q}ds\\
				&\leq C\int_{0}^{t}\mathbb{E}|e_{i}(s\wedge\tau_i)|^{p}ds+C\int_{0}^{T}[\mathbb{E}(\textbf{I}_{|\varpi_{s}^{i,N}|>\lambda^{-1}(\kappa(\varDelta t))}|\varpi_{s}^{i,N}|^{pq/(q-p\varepsilon)})]^{(q-p\varepsilon)/q}ds.\\
			\end{aligned}
		\end{equation}
		Then by applying the fundamental bridge and Chebyshev's inequality, we obtain
		\begin{equation}
			\label{j2211}
			\begin{aligned}
				J_{21}&\leq C\int_{0}^{t}\mathbb{E}|e_{i}(s\wedge\tau_i)|^{p}ds+C\int_{0}^{T}\Bigr([\mathbb{P}(|\varpi_{s}^{i,N}|>\lambda^{-1}(\kappa(\varDelta t)))]^{(q-p\varepsilon-p)/(q-p\varepsilon)}\\
				&\quad\times[\mathbb{E}|\varpi_{s}^{i,N}|]^{p/(q-p\varepsilon)}\Bigr)^{(q-p\varepsilon)/q}ds\\
				&\leq C\int_{0}^{t}\mathbb{E}|e_{i}(s\wedge\tau_i)|^{p}ds+C\int_{0}^{T}
				\Bigr(\frac{\mathbb{E}|\varpi_{s}^{i,N}|^{q}}{(\lambda^{-1}(\kappa(\varDelta t)))^{q}}\Bigr)^{(q-p\varepsilon-p)/q}ds\\
				&\leq  C\int_{0}^{t}\mathbb{E}|e_{i}(s\wedge\tau_i)|^{p}ds+C(\lambda^{-1}(\kappa(\varDelta t)))^{p\varepsilon+p-q}.
			\end{aligned}
		\end{equation}
		Then by applying Young inequality and Assumption 2.1, we have	
		\begin{equation}
			\label{j022}
			\begin{aligned}
				J_{22}&\leq C\int_{0}^{t}\mathbb{E}|e_{i}(s\wedge\tau_i)|^{p}ds
				+C\mathbb{E}\int_{0}^{t\wedge\tau_i}\Bigr(|e_{i}(s)|^{p-2}(1+|\varpi_{s}^{i,N}|^{2\varepsilon}+|\varpi_{[\eta s]}^{i,N}|^{2\varepsilon}\\
				&\quad+|{\omega^{*}}_s^{i,N}|^{2\varepsilon}+|{\tilde\omega}_s^{i,N}|^{2\varepsilon})
				\times(|\varpi_{s}^{i,N}-{\omega^{*}}_s^{i,N}|^{2}+|\varpi_{[\eta s]}^{i,N}-{\tilde\omega}_s^{i,N}|^{2})\\
				&\quad+\mathbb{W}_2^2(\upsilon_{s}^{\varpi,N},{\upsilon^{*}}_s^{\varpi,N})+\mathbb{W}_2^2(\upsilon_{[\eta s]}^{\varpi,N},{\tilde{\upsilon}}_s^{\varpi,N})\Bigr)ds\\
				&\leq C\int_{0}^{t}\mathbb{E}|e_{i}(s\wedge\tau_i)|^{p}ds+J_{221},
			\end{aligned}
		\end{equation}	
		where
		\begin{equation}
			\label{j22001}
			\begin{aligned}
				J_{221}&=C\mathbb{E}\int_{0}^{t\wedge\tau_i}|e_{i}(s)|^{p-2}(1+|\varpi_{s}^{i,N}|^{2\zeta}+|\varpi_{[\eta s]}^{i,N}|^{2\zeta}+|{\omega^{*}}_s^{i,N}|^{2\varepsilon}+|{\tilde\omega}_s^{i,N}|^{2\varepsilon})\\
				&\quad\times(|\varpi_{s}^{i,N}-{\omega^{*}}_s^{i,N}|^{2}+|\varpi_{[\eta s]}^{i,N}-{\tilde\omega}_s^{i,N}|^{2})ds.
			\end{aligned}
		\end{equation}		
		Then by applying Young inequality, H\"{o}lder inequality, Lemmas 4.2 and 4.3 and utilizing $pq/(q-p\varepsilon)\geq2$, we obtain
		\begin{equation}
			\label{j22010}
			\begin{aligned}
				J_{221}&\leq C\mathbb{E}\int_{0}^{t}|e_{i}(s\wedge\tau_i)|^{p}ds+C\mathbb{E}\int_{0}^{T}\Bigr(1+|\varpi_{s}^{i,N}|^{p\varepsilon}+|\varpi_{[\eta s]}^{i,N}|^{p\varepsilon}\\
				&\quad+|{\omega^{*}}_s^{i,N}|^{p\varepsilon}+|{\tilde\omega}_s^{i,N}|^{p\varepsilon}\Bigr)
				\times\Bigr(|\varpi_{s}^{i,N}-{\omega^{*}}_s^{i,N}|^{p}+|\varpi_{[\eta s]}^{i,N}-{\tilde\omega}_s^{i,N}|^{p}\Bigr)ds\\
				&\leq C\int_{0}^{t}\mathbb{E}|e_{i}(s\wedge\tau_i)|^{p}ds\\
				&\quad+C\int_{0}^{T}\Bigr[\mathbb{E}(1+|\varpi_{s}^{i,N}|^{q}+|\varpi_{[\eta s]}^{i,N}|^{q}+|{\omega^{*}}_s^{i,N}|^{q}+|{\tilde\omega}_s^{i,N}|^{q})\Bigr]^{p\varepsilon/q}\\
				&\quad\times\Bigr[\mathbb{E}|\varpi_{s}^{i,N}-{\omega^{*}}_s^{i,N}|^{pq/(q-p\varepsilon)}
				+\mathbb{E}|\varpi_{[\eta s]}^{i,N}-{\tilde\omega}_s^{i,N}|^{pq/(q-p\varepsilon)}\Bigr]^{(q-p\varepsilon)/q}ds\\
				&\leq C\int_{0}^{t}\mathbb{E}|e_{i}(s\wedge\tau_i)|^{p}ds+C\int_{0}^{T}
				\Bigr((\kappa(\varDelta t))^{pq/(q-p\varepsilon)}\Delta^{pq/2(q-p\varepsilon)}+\varDelta t\Bigr)^{(q-p\varepsilon)/q}ds\\
				&\leq C\int_{0}^{t}\mathbb{E}|e_{i}(s\wedge\tau_i)|^{p}ds+C((\kappa(\varDelta t))^{p}(\varDelta t)^{p/2}+(\varDelta t)^{(q-p\varepsilon)/q}).
			\end{aligned}
		\end{equation}
		Therefore by plugging in (\ref{j22010}) into (\ref{j022}) then substituting with (\ref{j2211}) and (\ref{j022}) into (\ref{j22}), we get
		\begin{equation}
			\label{second}
			\begin{aligned}
				J_{2}&\leq C\int_{0}^{t}\mathbb{E}|e_{i}(s\wedge\tau_i)|^{p}ds+C((\lambda^{-1}(\kappa(\varDelta t)))^{p\varepsilon+p-q}+(\kappa(\varDelta t))^{p}(\varDelta t)^{p/2}+(\varDelta t)^{(q-p\varepsilon)/q}).
			\end{aligned}
		\end{equation}
		By applying Assumption 2.1, Young inequality, Lemmas 4.2 and 4.3, we have	
		\begin{equation}
			\label{j33}
			\begin{aligned}
				J_3&=C\int_{0}^{t}\mathbb{E}|e_{i}(s\wedge\tau_i)|^{p}ds+C\int_{0}^{T}\mathbb{E}(|x_{s}^{i,N}-{\omega^{*}}_s^{i,N}|^{p}+|x_{[\eta s]}^{i,N}-{\tilde\omega}_s^{i,N})|^{p})ds\\
				&\leq C\mathbb{E}\int_{0}^{t}|e_{\Delta}(s\wedge\tau)|^{p}ds+C\int_{0}^{T}\mathbb{E}(|\varpi_{s}^{i,N}-{\omega^{*}}_s^{i,N}|^{p}+|\varpi_{[\eta s]}^{i,N}-{\tilde\omega}_s^{i,N}|^{p})ds\\
				&\leq C\mathbb{E}\int_{0}^{t}|e_{i}(s\wedge\tau_i)|^{p}ds+C((\kappa(\varDelta t))^{p}(\varDelta t)^{p/2}+\varDelta t).
			\end{aligned}	
		\end{equation}
		Then by plugging in (\ref{j11}), (\ref{second}) and (\ref{j33}) into (\ref{error4}), we get
		\begin{equation}
			\label{final1}
			\begin{aligned}
				\mathbb{E}|e_{i}(t)|^{p}&\leq C\int_{0}^{t}\mathbb{E}|e_{i}(s)|^{p}ds
				+C((\lambda^{-1}(\kappa(\varDelta t)))^{p\varepsilon+p-q}+(\kappa(\varDelta t))^{p}(\varDelta t)^{p/2}\\
				&\quad+(\varDelta t)^{(q-p\varepsilon)/q}).
			\end{aligned}
		\end{equation}
		Then Gronwall inequality implies
		\begin{equation}
			\label{final2}
			\mathbb{E}|e_{i}(t)|^{p}\leq C((\lambda^{-1}(\kappa(\Delta)))^{p\varepsilon+p-q}+(\kappa(\varDelta t))^{p}(\varDelta t)^{p/2}+(\varDelta t)^{(q-p\varepsilon)/q}).
		\end{equation}
		By combining (\ref{final2}) with (\ref{poc}), the required assertion (\ref{main66}) is attained. Furthermore, by utilizing (\ref{qq}), we get $\lambda^{-1}(r)=(r/C)^{1/(1+\varepsilon)}$. By substituting with $\lambda^{-1}(r)$ and (\ref{qqq}) into (\ref{main66}), we get (\ref{main660}) and the proof is complete.}}\\\\
\textbf{{\small Corollary 4.2.}} {\small\textit{ Let assumptions 2.1-2.3 hold and let $\lambda(r)$ and $\kappa(\varDelta t)$ be defined in (\ref{qq}) and (\ref{qqq}). Then for any 
		\begin{equation}
			\label{rrr}
			p\in[2,q_0),\quad q\in((1+\varepsilon)p\vee q_0,q^*)\quad\text{and}\quad \gamma\in (0,1/4\wedge 1/q],
		\end{equation}		
		we have
		\begin{equation}
			\label{main6600}
			\begin{aligned}
				\sup _{t\in[0, T]}\mathbb{E}|x_{t}^{i}-\varpi_{t}^{i,N}|^{p}
				\leq C
				\left\{\begin{array}{ll}
					N^{-1 / 2}+D^*, & \text { if } p>d/2, \\
					N^{-1 / 2} \log(1+N)+D^*, & \text { if } p=d/2,\\
					N^{-p / d}+D^*, & \text { if }p\in[2,d/2),
				\end{array}\right.
			\end{aligned}
		\end{equation}
		where
		\begin{equation}
			\label{main200}
			D^*=C\varDelta t^{[\gamma(q-(1+\varepsilon)p)/(1+\varepsilon)]\wedge[(q-\varepsilon p)/q]}
		\end{equation}
		for all $\varDelta t\in(0,\varDelta t^*]$.}}\\\\
\textbf{{\small Proof.}} {\small  Firstly, replacing condition $q>(1+\varepsilon)q_0$ by a weaker one $q>(1+\varepsilon)p\vee q_0$ does not affect the results in Theorem 4.1 , but it makes some relaxation and flexibility in the choice of $q$. Therefore, by utilizing (\ref{rrr}), we have 
	\begin{equation*}
		(\gamma(q-(1+\varepsilon)p)/(1+\varepsilon))<p(1-2\gamma)/2.
	\end{equation*}	
	Then, by (\ref{main660}), we get (\ref{main6600}). The proof is complete.}\\\\

\textbf{{\small Remark 4.1.}} {\small\textit{ Let we assume that $q^*$ is allowed to be large. Then by Corollary (4.2), when $q$ is sufficiently large relative to $\varepsilon p$, the term $D^*$ will be dominated by $C(\varDelta t)^{(\gamma(q-(1+\varepsilon)p)/(1+\varepsilon))}$, where $\gamma(q-(1+\varepsilon)p)/(1+\varepsilon)$ is proportional to $\varepsilon$. Therefore, by substituting $\gamma=1/4\wedge 1/q$ into $D^*$, we get
		\begin{equation}
			D^*=C(\varDelta t)^{((1/(1+\varepsilon))-(p/q))}.
		\end{equation}			
		Then, by letting $q\rightarrow\infty$, we get 
		\begin{equation}
			\label{main66000}
			\begin{aligned}
				\sup _{t\in[0, T]}\mathbb{E}|x_{t}^{i}-\varpi_{t}^{i,N}|^{p}
				\leq C
				\left\{\begin{array}{ll}
					N^{-1 / 2}+C(\varDelta t)^{1/(1+\varepsilon)}, & \text { if } p>d/2, \\
					N^{-1 / 2} \log(1+N)+C(\varDelta t)^{1/(1+\varepsilon)}, & \text { if } p=d/2,\\
					N^{-p / d}+C(\varDelta t)^{1/(1+\varepsilon)}, & \text { if }p\in[2,d/2),
				\end{array}\right.
			\end{aligned}
		\end{equation}
		for all $\varDelta t\in(0,\varDelta t^*]$ and this can be considered as the optimal convergence rate of the truncated Euler-Maruyama scheme in the case of jumps. Furthermore, it should be mentioned that this is significantly different from the case without jumps where the results will be 
		\begin{equation}
			\label{main660000}
			\begin{aligned}
				\sup _{t\in[0, T]}\mathbb{E}|x_{t}^{i}-\varpi_{t}^{i,N}|^{p}
				\leq C
				\left\{\begin{array}{ll}
					N^{-1 / 2}+C(\varDelta t)^{p(1-2\gamma)/2}, & \text { if } p>d/2, \\
					N^{-1 / 2} \log(1+N)+C(\varDelta t)^{p(1-2\gamma)/2}, & \text { if } p=d/2,\\
					N^{-p / d}+C(\varDelta t)^{p(1-2\gamma)/2}, & \text { if }p\in[2,d/2),
				\end{array}\right.
			\end{aligned}
		\end{equation}
		for all $\varDelta t\in(0,\varDelta t^*]$ and $\gamma\in(0,1/4]$ where the $\mathcal{L}^{p}(p\geq2)$ convergence rate is close to $p/2$ when there is no jumps in the stochastic model \cite{guo2018note}.}}

\section{{\normalsize Stability analysis}}
{\small
	In this section, the almost sure exponential stability of the exact and numerical solutions of our stochastic model will be addressed.\\\\
	\textbf{Assumption 5.1.} \textit{ There exist positive constants $c_i$, $i=1,2,3,4$ such that
		\begin{equation}
			\label{5.1}
			\begin{aligned}
				2\langle \vartheta,f(\vartheta,\xi,\nu,\nu^*)\rangle&+\|g(\vartheta,\xi,\nu,\nu^*)\|^{2}+\int_{Z}|h(\vartheta,\xi,\nu,\nu^*,z)|^{2}\pi(dz)\\
				&\quad\leq-c_1|\vartheta|^2+c_2|\xi|^2+c_3\mathcal{W}_2^2(\nu,\delta_0)+c_4	\mathcal{W}_2^2(\nu^*,\delta_0)
			\end{aligned}	
		\end{equation}
		for any $\vartheta,\xi\in\mathbb{R}^{d}$ and $\nu,\nu^*\in\mathcal{P}_2(\mathbb{R}^d)$}.\\\\
	\textbf{Assumption 5.2.} \textit{ There exist positive constants $b_i $, $i=1,2,3,4$ such that
		\begin{equation}
			\label{5.2}
			\begin{aligned}
				|f(\vartheta,\xi,\nu,\nu^*)|^2\leq b_1|\vartheta|^2+b_2|\xi|^2+b_3\mathcal{W}_2^2(\nu,\delta_0)+b_4	\mathcal{W}_2^2(\nu^*,\delta_0)
			\end{aligned}	
		\end{equation}
		for any $\vartheta,\xi\in\mathbb{R}^{d}$ and $\nu,\nu^*\in\mathcal{P}_2(\mathbb{R}^d)$}.\\\\
	By utilizing Assumptions 5.1 and 5.2, we have 
	\begin{equation}
		\label{5.11}
		\begin{aligned}
			2\langle \vartheta,f_\Delta(\vartheta,\xi,\nu,\nu^*)\rangle&+\|g_\Delta(\vartheta,\xi,\nu,\nu^*)\|^{2}+\int_{Z}|h(\vartheta,\xi,\nu,\nu^*,z)|^{2}\pi(dz)\\
			&\quad\leq-c_1|\vartheta|^2+c_2|\xi|^2+c_3\mathcal{W}_2^2(\nu,\delta_0)+c_4	\mathcal{W}_2^2(\nu^*,\delta_0),
		\end{aligned}	
	\end{equation}
	and
	\begin{equation}
		\label{5.22}
		\begin{aligned}
			|f_\Delta(\vartheta,\xi,\nu,\nu^*)|^2\leq b_1|\vartheta|^2+b_2|\xi|^2+b_3\mathcal{W}_2^2(\nu,\delta_0)+b_4	\mathcal{W}_2^2(\nu^*,\delta_0).
		\end{aligned}	
	\end{equation}
	for any $\vartheta,\xi\in\mathbb{R}^{d}$ and $\nu,\nu^*\in\mathcal{P}_2(\mathbb{R}^d)$}.\\\\
\textbf{Definition 5.1} \cite{liu2023stability}. The solution $x_t^{i,N}$ to Eq.(\ref{modelmv2}) possesses the almost sure exponential stability property if there is a positive constant $\gamma_1$ such that
\begin{equation}
	\limsup_{t\rightarrow \infty}\frac{1}{t}\log(\frac{1}{N}\sum_{i=1}^{N}|x_t^{i,N}|)\leq -\gamma_1 \quad, a.s.
\end{equation}	
\textbf{Definition 5.2} \cite{liu2023stability}. The numerical solution $ \{	\varpi_{l}^{i,N}\}_{l\geq1} $ to Eq.(\ref{TEM}) possesses the almost sure exponential stability property if there is a positive constant $\gamma_2$ such that,
\begin{equation}
	\limsup_{n\rightarrow \infty}\frac{1}{l\varDelta t}\log(\frac{1}{N}\sum_{i=1}^{N}|	\varpi_{l}^{i,N}|)\leq -\gamma_2 \quad a.s.
\end{equation}	
\textbf{Theorem 5.1.} \textit{Let Assumptions 2.1 and 5.1 hold and $(c_1-c_3)-(c_2-c_4)\geq0$. Then the solution $x_t^{i,N} $ to Eq.(\ref{modelmv2}) possesses the property of almost sure exponential stability
	\begin{equation}
		\limsup_{t\rightarrow \infty}\frac{1}{t}\log(\frac{1}{N}\sum_{i=1}^{N}|x_t^{i,N}|)\leq -\frac{\gamma^*}{2}\quad, a.s.,
	\end{equation}	
	where $\gamma^*=\min(1,(c_1-c_3)-(c_2-c_4))$.}\\\\
\textbf{Proof.} The proof of the theorem can be attained by the same concept and techniques discussed in \cite{liao1997almost}.\\\\
\textbf{Theorem 5.2.} \textit{Let Assumptions 2.1, 5.1 and 5.2 hold and $(c_1-c_3)>(c_2+c_4)(1+[1/\eta])$. Then, there exists a constant $ \varDelta t^{*} $ such that for all $ \varDelta t \in (0,\varDelta t^*) $ and initial data $\varpi_0^{i,N} $, the approximate solution defined by Eq.(\ref{TEM}) possesses the property 
	\begin{equation}
		\label{th1stability}
		\lim_{l\rightarrow\infty}\frac{1}{l\varDelta t}\log	\big(\frac{1}{N}\sum_{i=1}^{N}|\varpi_{l}^{i,N}|^2\big)\leq -\gamma\quad a.s.
	\end{equation}	
	where $ \gamma=\min(1,\varsigma^{*} ) $ with $ \varsigma^{*} $ satisfying
	\begin{equation}
		\begin{aligned}
			\varsigma^{*}+(b_1+b_3)\varDelta t)+c_3-c1+((b_2+b_4)\varDelta t+c_2+c_4)(1+[1/\eta])=0
		\end{aligned}
	\end{equation}
	and	
	\begin{equation}
		\lim_{\varDelta t\rightarrow 0}\varsigma^{*}=c_1-c_3-(c_2+c_4)(1+[1/\eta]).
\end{equation}}
\textbf{Proof.} From (\ref{TEM}), we can obtain
\begin{equation}
	\label{scheme110}
	\begin{aligned}		
		|\varpi_{l+1}^{j,N}|^{2}&=|\varpi_{l}^{j,N}|^{2}
		+\varDelta(2\langle\varpi_{l}^{j,N},f_{\Delta}(\varpi_{l}^{j,N},\varpi_{[\eta l]}^{j,N},\upsilon_{l}^{\varpi,N},\upsilon_{[\eta l]}^{\varpi,N})\rangle \\
		&\quad+\|g_{\Delta}(\varpi_{l}^{j,N},\varpi_{[\eta l]}^{i,N},\upsilon_{l}^{\varpi,N},\upsilon_{[\eta l]}^{\varpi,N})\|^{2}\\
		&\quad+\int_Z |h(\varpi_{l}^{j,N},\varpi_{[\eta l]}^{j,N},\upsilon_{l}^{\varpi,N},\upsilon_{[\eta l]}^{\varpi,N},z)|^{2}\pi(dz)\\
		&\quad+|f_{\Delta}(\varpi_{l}^{j,N},\varpi_{[\eta l]}^{j,N},\upsilon_{l}^{\varpi,N},\upsilon_{[\eta l]}^{\varpi,N})|^2\varDelta t)+M^j_{l},
	\end{aligned}		
\end{equation}
where
\begin{equation}
	\begin{aligned}
		M_{l}^j&=\|g_{\Delta}(\varpi_{l}^{j,N},\varpi_{[\eta l]}^{j,N},\upsilon_{l}^{\varpi,N},\upsilon_{[\eta l]}^{\varpi,N})\|^{2}(|\varDelta W_l^{j}|^2-\varDelta t)\\
		&\quad+\int_{Z}|h(\varpi_{l}^{j,N},\varpi_{[\eta l]}^{j,N},\upsilon_{l}^{\varpi,N},\upsilon_{[\eta l]}^{\varpi,N},z)|^{2}(|\widetilde{N}^j(\varDelta t,dz)|^{2}-\pi(dz)\varDelta t)\\
		&\quad2\langle \varpi_{l}^{j,N},g_{\Delta}(\varpi_{l}^{j,N},\varpi_{[\eta l]}^{j,N},\upsilon_{l}^{\varpi,N},\upsilon_{[\eta l]}^{\varpi,N})\varDelta W_l^j\rangle\\
		&\quad+2\langle f_{\Delta}(\varpi_{l}^{j,N},\varpi_{[\eta l]}^{j,N},\upsilon_{l}^{\varpi,N},\upsilon_{[\eta l]}^{\varpi,N})\varDelta t,g_{\Delta}(\varpi_{l}^{j,N},\varpi_{[\eta l]}^{j,N},\upsilon_{l}^{\varpi,N},\upsilon_{[\eta l]}^{\varpi,N})\varDelta W_l^j\rangle\\
		&\quad2\langle \varpi_{l}^{j,N},\int_{Z}h(\varpi_{l}^{j,N},\varpi_{[\eta l]}^{j,N},\upsilon_{l}^{\varpi,N},\upsilon_{[\eta l]}^{\varpi,N},z)\widetilde{N}^j(\varDelta t,dz)\rangle\\
		&\quad+2\langle f_{\Delta}(\varpi_{l}^{j,N},\varpi_{[\eta l]}^{j,N},\upsilon_{l}^{\varpi,N},\upsilon_{[\eta l]}^{\varpi,N})\varDelta t,\\
		&\quad\int_{Z}h(\varpi_{l}^{j,N},\varpi_{[\eta l]}^{j,N},\upsilon_{l}^{\varpi,N},\upsilon_{[\eta l]}^{\varpi,N},z)\widetilde{N}^j(\varDelta t,dz)\rangle\\
		&\quad+2\langle g_{\Delta}(\varpi_{l}^{j,N},\varpi_{[\eta l]}^{j,N},\upsilon_{l}^{\varpi,N},\upsilon_{[\eta l]}^{\varpi,N})\varDelta W_l^j,\\
		&\quad\int_{Z}h(\varpi_{l}^{j,N},\varpi_{[\eta l]}^{j,N},\upsilon_{l}^{\varpi,N},\upsilon_{[\eta l]}^{\varpi,N},z)\widetilde{N}^j(\varDelta t,dz)\rangle.
	\end{aligned}	
\end{equation}	
Then by applying Assumptions 5.1 and 5.2, we get 
\begin{equation}
	\label{equality2}
	\begin{aligned}
		|\varpi_{l+1}^{j,N}|^{2}&\leq|\varpi_{l}^{j,N}|^2+(b_1\varDelta t-c_1)\varDelta t|\varpi_{l}^{j,N}|^{2}
		+(b_3\varDelta t+c_3)\varDelta t\frac{1}{N}\sum_{i=1}^{N}|\varpi_{l}^{i,N}|^2\\
		&\quad+(b_2\varDelta t+c2)\varDelta t|\varpi^{j,N}_{[\eta l]}|^{2}
		+(b_4\varDelta t+c_4)\varDelta t\frac{1}{N}\sum_{i=1}^{N}|\varpi_{[\eta l]}^{i,N}|^2
		+M_l^j.
	\end{aligned}	
\end{equation}	
Then for any constant $\varsigma>0$, we have the following
\begin{equation}
	\label{inequality}
	\begin{aligned}
		e^{\varsigma(l+1)\varDelta t}|\varpi_{l+1}^{j,N}|^2-e^{\varsigma l\varDelta t}|\varpi_{l}^{j,N}|^{2}&=e^{\varsigma(l+1)\varDelta t}(|\varpi_{l+1}^{j,N}|^2-|\varpi_{l}^{j,N}|^2)\\
		&\quad+(1-e^{-\varsigma\varDelta t})e^{\varsigma(l+1)\varDelta t}|\varpi_{l}^{j,N}|^{2}.
	\end{aligned}	
\end{equation}	
By plugging in (\ref{equality2}) into (\ref{inequality}), we obtain
\begin{equation}
	\label{main40}
	\begin{aligned}
		e^{\varsigma(l+1)\varDelta t}&|\varpi_{l+1}^{j,N}|^2-e^{\varsigma l\varDelta t}|\varpi_{l+1}^{j,N}|^{2}\\
		&\leq e^{\varsigma(l+1)\varDelta t}\Biggl[E(\varsigma,\varDelta t)\varDelta t|\varpi_{l}^{j,N}|^{2}+F(\varsigma,\varDelta t)\varDelta t |\varpi^{j,N}_{[\eta l]}|^{2}\\
		&\quad+G(\varsigma,\varDelta t)\varDelta t\frac{1}{N}\sum_{i=1}^{N}|\varpi_{l}^{i,N}|^2+H(\varsigma,\varDelta t)\varDelta t \frac{1}{N}\sum_{i=1}^{N}|\varpi_{[\eta l]}^{i,N}|^2+M^j_{l}\Biggl],
	\end{aligned}	
\end{equation}
where
\begin{equation}
	E(\varsigma,\varDelta t)=(\varsigma+b_1\varDelta t-c_1),\quad F(\varsigma,\varDelta t)=(b_2\varDelta t+c_2),
\end{equation}
and
\begin{equation}
	G(\varsigma,\varDelta t)=(b_3\varDelta t+c_3),\quad H(\varsigma,\varDelta t)=(b_4\varDelta t+c_4).
\end{equation}		
Taking the summation from $ k=0 $ to $ k=l-1 $ to both sides of (\ref{main40}) results in
\begin{equation}
	\label{main400}
	\begin{aligned}
		e^{\varsigma l\varDelta t}|\varpi_{l}^{j,N}|^2&\leq |\varpi^{j,N}_0|^{2}+E(\zeta,\varDelta t)\varDelta t \sum_{k=0}^{l-1}e^{\varsigma(k+1)\varDelta t}|\varpi_{k}^{j,N}|^{2}\\
		&\quad+F(\varsigma,\varDelta t)\varDelta t \sum_{k=0}^{l-1}e^{\varsigma(k+1)\varDelta t} |\varpi^{j,N}_{[\eta k]}|^{2}\\
		&\quad+G(\varsigma,\varDelta t)\varDelta t\frac{1}{N}\sum_{k=0}^{l-1}\sum_{i=1}^{N}e^{\varsigma(k+1)\varDelta t}|\varpi_{k}^{i,N}|^2\\
		&\quad+H(\varsigma,\varDelta t)\varDelta t \frac{1}{N}\sum_{k=0}^{l-1}\sum_{i=1}^{N}e^{\varsigma(k+1)\varDelta t}|\varpi_{[\eta k]}^{i,N}|^2\\
		&\quad+\sum_{k=0}^{l-1}e^{\varsigma(k+1)\varDelta t} M^j_{k} 
	\end{aligned}	
\end{equation}
where $\sum_{k=0}^{l-1}e^{\varsigma(k+1)\varDelta t} M^j_{k} $ is a martingale. Furthermore it can be concluded that 
\begin{equation}
	\label{main421}
	\begin{aligned}
		\sum_{k=0}^{l-1}&e^{\varsigma(k+1)\varDelta t} |\varpi^{j,N}_{[\eta k]}|^{2}\\
		&\leq (1+[1/\eta])\Biggl(\sum_{k=0}^{l-1}e^{\varsigma(k+1)\varDelta t}|\varpi_{k}^{j,N}|^{2}
		-\sum_{k=[\eta(l-1)]+1}^{l-1}e^{\varsigma(k+1)\varDelta t}|\varpi^{j,N}_{ k}|^{2}\Biggr),
	\end{aligned}	
\end{equation}	
and
\begin{equation}
	\label{main422}
	\begin{aligned}
		\sum_{k=0}^{l-1}&\sum_{i=1}^{N}e^{\varsigma(k+1)\varDelta t}|\varpi_{[\eta k]}^{i,N}|^2\\
		&\leq (1+[1/\eta])\Biggl(\sum_{k=0}^{l-1}\sum_{i=1}^{N}e^{\varsigma(k+1)\varDelta t}|\varpi_{k}^{i,N}|^{2}
		-\sum_{k=[\eta(l-1)]+1}^{l-1}\sum_{i=1}^{N}e^{\varsigma(k+1)\varDelta t}|\varpi^{i,N}_{k}|^{2}\Biggr).
	\end{aligned}	
\end{equation}
By plugging in (\ref{main421}) and (\ref{main422}) into (\ref{main400}), we obtain
\begin{equation}
	\label{help2}
	\begin{aligned}
		\sum_{i=1}^{N}e^{\varsigma l\varDelta t}|\varpi^{i,N}_l|^2&\leq \sum_{i=1}^{N}|\varpi^{i,N}_0|^{2}+T(\varsigma,\varDelta t)\varDelta t \sum_{k=0}^{l-1}\sum_{i=1}^{N}e^{\varsigma(k+1)\varDelta t}|\varpi^{i,N}_k|^{2}\\
		&\quad-F(\varsigma,\varDelta t)\varDelta t(1+[1/\eta]) \sum_{k=[\eta(l-1)]+1}^{l-1}\sum_{i=1}^{N}e^{\zeta(k+1)\varDelta t}|\varpi^{i,N}_k|^{2}\\
		&\quad-H(\varsigma,\varDelta t)\varDelta t(1+[1/\eta]) \sum_{k=[\eta(l-1)]+1}^{l-1}\sum_{i=1}^{N}e^{\zeta(k+1)\varDelta t} |\varpi^{i,N}_k|^{2}\\
		&\quad+\sum_{k=0}^{l-1}\sum_{i=1}^{N}e^{\varsigma(k+1)\varDelta t} M^i_{k},
	\end{aligned}	
\end{equation}
where $ T(\zeta,\varDelta t):=E(\varsigma,\varDelta t)+G(\varsigma,\varDelta t)+(F(\varsigma,\varDelta t)+H(\varsigma,\varDelta t))(1+[1/\eta])$. Let
\begin{equation}
	\varDelta t^{*}=\frac{(c_1-c_3)-(c_2+c_4)(1+[1/\eta])}{(b_1+b_3)+(b_2+b_4)(1+[1/\eta])}.
\end{equation}	
Then for any $ \varDelta t < \varDelta t^{*} $, 
\begin{equation}
	T(0,\varDelta t)=(b_1+b_3)\varDelta t+(c_3-c_1)+((b_2+b_4)\varDelta t+c_2+c_4) (1+[1/\eta])<0.	
\end{equation}
It is also noted that $ \frac{d}{d\zeta}T(\varsigma,\varDelta t)>0 $. So, there exists a unique constant $ \varsigma^{*}>0 $ satisfying $ T(\varsigma^{*},\varDelta t)=0 $. Also, we have
\begin{equation}
	\lim_{\varDelta t\rightarrow 0}T((c_1-c_3)-(c_2+c_4)(1+[1/\eta]),\varDelta t)=0.
\end{equation}	 
Therefore, we conclude that 
\begin{equation}
	\lim_{\varDelta t\rightarrow 0}\varsigma^{*}=(c_1-c_3)-(c_2+c_4)(1+[1/\eta]).
\end{equation}
Using the definition of $ \gamma $, if $ \varsigma^{*}>1 $, then $ \gamma=1 $ and $ T(\gamma,\varDelta t)<0 $. However, if $ \varsigma^{*}<1$ then $ \gamma=\zeta^{*} $ and $ T(\gamma,\varDelta t)=0.$ By applying the discrete semi-martingale convergence theorem \cite{mao2006stochastic}, we have for any initial data $ \varpi_0^{i,N} $
\begin{equation}
	\limsup_{l\rightarrow \infty}\Biggl[\sum_{i=1}^{N}|\varpi_0^{i,N}|^{2}+\sum_{k=0}^{l-1}\sum_{i=1}^{N}e^{\varsigma(k+1)\varDelta t} M^i_{k}\Biggr]<\infty,
\end{equation}
which yields
\begin{equation}
	\begin{aligned}
		\limsup_{l\rightarrow \infty}e^{\gamma l\varDelta t}\frac{1}{N}\sum_{i=1}^{N}|\varpi_l^{i,N}|^{2}\leq\infty,\quad a.s.
	\end{aligned}
\end{equation}	
Therefore there is a finite positive random variable $ \epsilon $ such that
\begin{equation}
	\limsup_{l\rightarrow \infty}e^{\gamma l\varDelta t}\frac{1}{N}\sum_{i=1}^{N}|\varpi_l^{i,N}|^{2}<\epsilon,
\end{equation}
which implies the required assertion (\ref{th1stability}). The proof is complete.
\section{{\normalsize Numerical examples}}
{\small
\textbf{{\small Example 6.1.}} Consider the following stochastic proportional delay Mckean-Vlasov model with L\'evy jumps
\begin{equation}
	\label{example}
	\begin{aligned}
		dx_t=(-2x_t^{5}+\mathbb{E}[x_t])dt+(x_{t}^{2}+x_{\eta t}\sin^{2}(x_t))dW(t)
		+\int_{Z}x_tz\widetilde{N}(dt,dz),\quad t\geq0
	\end{aligned}
\end{equation}	
with initial data $x_0\sim N(0,1)$ where $N(0,1)$ is the standard normal distribution, $\eta=0.4$ and compensator is given by $\pi(dz)dt=f(z)dzdt$ where $f(z)\sim N(0,1)$. It can be noticed the coefficients of Eq.(\ref{example}) satisfy Assumptions 2.1-2.3 and 
\begin{equation}
	\sup_{|\vartheta|\vee|\xi|\leq r}\left(|f(\vartheta,\xi,\nu,\nu^*)|\vee\left\|g(\vartheta,\xi,\nu,\nu^*)\right\|\right) \leq 2r^5,\quad r\geq1.
\end{equation}	
Therefore, we can select $\lambda(r)=r^5$, $r\geq0$ and $\kappa(\varDelta t)=\varDelta t^{-1/10}$ for the objection of truncation and the numerical algorithm Eq.(\ref{TEM}) can be applied to our example Eq.(\ref{example}). Moreover, to approximate the law of $x_t$ at each time step by its empirical distribution, we apply the particle method with number of particles $N=100$. Because Eq.(\ref{example}) has no exact solution, it can be pretended that the output of our numerical algorithm Eq.(\ref{TEM}) with very small step-size $dt=2^{-13}$ is the exact solution of Eq.(\ref{example}) for the error comparison with other step sizes $\varDelta t=2^{i}dt$, $1\leq i\leq 4$. Furthermore sample paths of Brownian motion are generated with step size $2^{-13}$ and Monte Carlo approach with 1000 sample paths is utilized to find the root-mean-square error (RMSE) as follows
\begin{equation}
	\text{RMSE}=\sqrt{\frac{1}{MN}\sum_{j=1}^{M}\sum_{i=1}^{N}|\varpi_T^{i,N,*,j}-\varpi_T^{i,N,j}|^2},
\end{equation}	
where $\varpi_T^{i,N,*}$ and $\varpi_T^{i,N}$ represents the exact and approximated solutions, respectively, for the $i-$th particle at terminal time $T=1$. Figure \ref{fig:graph1} shows the plot of the RMSE against $\varDelta t$ on a log-log scale and as a referenced dashed line of slope $1/2$ is also added. The convergence rate is approximately one half.\\
\begin{figure}[!h]
	\centering
	\includegraphics[width=.5\textwidth]{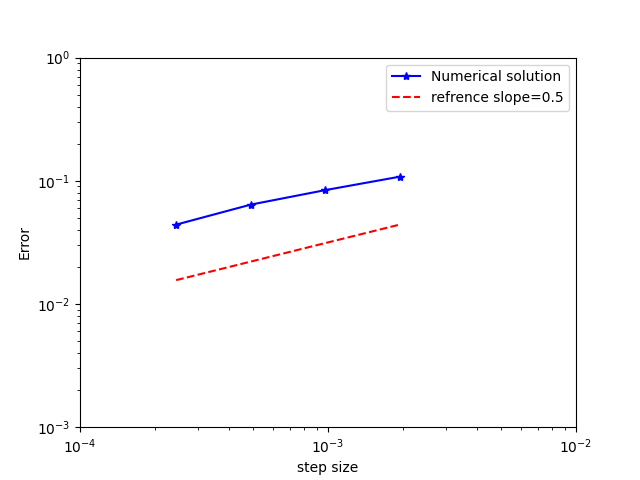}
	\caption{{\small Log-log plot of RMSE versus $\varDelta t $ for Eq.(\ref{example}).}}
	\label{fig:graph1}
\end{figure}\\\\

\textbf{{\small Example 6.2.}} Consider the following stochastic proportional delay Mckean-Vlasov model with L\'evy jumps
\begin{equation}
	\label{example6.2}
	\begin{aligned}
		dx_t=(-2x_t^{4}+\mathbb{E}[x_t^2])dt+(x_{t}^{5}+x_{\eta t}\cos^{2}(x_t))dW(t)
		+\int_{Z}(1+x_t)z\widetilde{N}(dt,dz),\quad t\geq0
	\end{aligned}
\end{equation}	
with initial data $x_0\sim N(0,1)$ where $N(0,1)$ is the standard normal distribution, $\eta=0.8$ and compensator is given by $\pi(dz)dt=f(z)dzdt$ where $f(z)\sim N(0,1)$. It can be noticed the coefficients of Eq.(\ref{example6.2}) satisfy Assumptions 2.1-2.3 and 
\begin{equation}
	\sup_{|\vartheta|\vee|\xi|\leq r}\left(|f(\vartheta,\xi,\nu,\nu^*)|\vee\left\|g(\vartheta,\xi,\nu,\nu^*)\right\|\right) \leq 2r^6,\quad r\geq1.
\end{equation}	
Therefore, we can select $\lambda(r)=r^6$, $r\geq0$ and $\kappa(\varDelta t)=\varDelta t^{-1/10}$ for the objection of truncation. Therefore, the truncated EM algorithm given by Eq.(\ref{TEM}) can be applied to our example Eq.(\ref{example6.2}). Moreover, we will present a comparison with the tamed EM method \cite{gao2024convergence} which is defined by	
\begin{equation}
	\label{TEM0}
	\begin{aligned}	
		Y_{l+1}^{i,N}&=Y_{l}^{i,N}+ \Delta \frac{f(Y_{l}^{i,N},Y_{[\eta l]}^{i,N},\upsilon_{l}^{Y,N},\upsilon_{[\eta l]}^{Y,N})}{1+\Delta^\beta|f_{\Delta}(Y_{l}^{i,N},Y_{[\eta l]}^{i,N},\upsilon_{l}^{Y,N},\upsilon_{[\eta l]}^{Y,N})|}
		+g(Y_{l}^{i,N},Y_{[\eta l]}^{i,N},\upsilon_{l}^{Y,N},\upsilon_{[\eta l]}^{Y,N})\Delta W_{l}^{i}\\
		&\quad+\int_{Z}h(Y_{l}^{i,N},Y_{[\eta l]}^{i,N},\upsilon_{l}^{Y,N},\upsilon_{[\eta l]}^{Y,N},z)\widetilde{N}^i(\varDelta t,dz),
	\end{aligned}
\end{equation}
with initial value $Y_{0}^{i,N}=x^{i,N}_0$ and $\beta=1/2$ where $\Delta W_{l}^{i}=W_{t_{l+1}}^{i}-W_{t_{l}}^{i}$. Also,
\begin{equation}
	\upsilon_{l}^{Y,N}(\cdot):=\frac{1}{N}\sum_{i=1}^N \delta_{Y_{l}^{i,N}}(\cdot)\quad\text{and}\quad
	\upsilon_{[\eta l]}^{Y,N}(\cdot):=\frac{1}{N}\sum_{i=1}^N \delta_{Y_{[\eta l]}^{i,N}}(\cdot).
\end{equation}
Then, we consider numerical approximations generated by Eq.(\ref{TEM}) and Eq.(\ref{TEM0}) with very small step-size $dt=2^{-12}$ as the exact solutions of Eq.(\ref{example6.2}) for the error comparison with other step sizes $\varDelta t\in\{2^{-11},2^{-10},2^{-9},2^{-8}\}$. We also apply the particle method with $N=500$ to approximate the law of $x_t$ at each time step. Furthermore sample paths of Brownian motion are generated with step size $2^{-12}$ and Monte Carlo approach with $ M=1000$ sample paths is utilized to find the root-mean-square error (RMSE) as follows
\begin{equation}
	\text{RMSE}=\sqrt{\frac{1}{MN}\sum_{j=1}^{M}\sum_{i=1}^{N}|Q_T^{i,N,*,j}-Q_T^{i,N,j}|^2},
\end{equation}	
where $Q_T^{i,N,*}$ and $Q_T^{i,N}$ represents the exact and approximated solutions, respectively, generated by numerical approximations Eq.(\ref{TEM}) and Eq.(\ref{TEM0}) for the $i-$th particle at terminal time $T=4$. Table 1 lists the root mean square errors of these approximations at $T=4$. From Table \ref{tab:data}, it can noticed that the truncated Euler-Maruyama method has alittle bit higher accuracy and less mean square error than the tamed EM scheme. Also, when the step size decreases, the root mean square error in both numerical approximations decreases.\\\\
\begin{table}[!h]
	\centering
	\caption[loftitle]{{RMSE of the numerical approximations}}
	\label{table1}\label{tab:data}
	\begin{tabular}{ccccccc}
		\hline
		$\varDelta t$&&& Tamed EM &&& Truncated EM\\
		$2^{-8}$&&&0.1761 &&& 0.1522 \\
		$2^{-9}$&&&0.1421 & &&0.0763 \\
		$2^{-10}$&&&0.0928 &&& 0.0541 \\
		$2^{-11}$&&&0.0654 &&& 0.0324 \\
		\hline
	\end{tabular}	
\end{table}

\textbf{{\small Example 6.3.}} Consider the following stochastic proportional delay Mckean-Vlasov model with L\'evy jumps
\begin{equation}
	\label{exampleee}
	\begin{aligned}
		dx_t=(-x_t^{7}+\mathbb{E}[x_t])dt+x_{\eta t}^{2}dW(t)
		+\int_{Z}(x_t+x_{\eta t})z\widetilde{N}(dt,dz),\quad t\geq0
	\end{aligned}
\end{equation}		
with initial data $x_0\sim N(0,1)$ where $N(0,1)$ is the standard normal distribution, $\eta=0.2$ and compensator is given by $\pi(dz)dt=f(z)dzdt$ where $f(z)\sim N(0,1)$. The coefficients of Eq.(\ref{exampleee}) satisfy Assumptions 2.1-2.3. By letting $p=2$ and $q_0=5$ and also selecting $\lambda(r)=r^7$, $r\geq0$, $\varepsilon=4$ and $\kappa(\varDelta t)=\varDelta t^{-1/50}$ and $q=50$, then by Corollary (4.2), the numerical solution will convergence to the true solution in $\mathcal{L}^{2}$ with convergence rate $1/(1+\varepsilon)-p/q=0.16$. As we see, Eq.(\ref{exampleee}) does not have exact solution, so we apply the truncated EM method with step size $dt=2^{-16}$ and consider its output as the exact solution for Eq.(\ref{exampleee}) for the error comparison with other step sizes $\varDelta t=\{2^{-12},2^{-13},2^{-14},2^{-15}\}$. We also apply the particle method with number of particles $N=400$ to approximate the distribution of $x_t$. Sample paths  of Brownian motion are generated with step size $2^{-17}$ and Monte Carlo approach with $M=1000$ sample paths is utilized to find the root-mean-square error (RMSE) as follows	
\begin{equation}
	\text{RMSE}=\sqrt{\frac{1}{MN}\sum_{j=1}^{M}\sum_{i=1}^{N}|\varpi_T^{i,N,*,j}-\varpi_T^{i,N,j}|^2},
\end{equation}	
where $\varpi_T^{i,N,*}$ and $\varpi_T^{i,N}$ represents the exact and approximated solutions, respectively, for the $i-$th particle at terminal time $T=1$. Figure \ref{fig:graph3} shows the plot of the RMSE against $\varDelta t$ on a log-log scale.\\\\
\begin{figure}[!h]
	\centering
	\includegraphics[width=.5\textwidth]{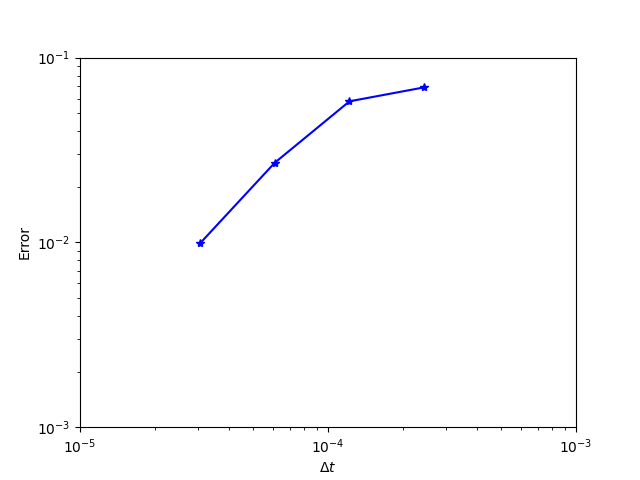}
	\caption{{\small Log-log plot of RMSE versus $\varDelta t $ for Eq.(\ref{exampleee}).}}
	\label{fig:graph3}
\end{figure}\\\\			
\textbf{{\small Example 6.4.}} Consider the following stochastic proportional delay Mckean-Vlasov model with L\'evy jumps
\begin{equation}
	\label{example02}
	\begin{aligned}
		dx_t&=(-x^{3}_t-x^{5}_t)dt+(x_t\cos^{2}(x_{\eta t})-0.5\int_{\mathbb{R}^d}y\nu(dy))dW(t)\\
		&\quad+\int_{Z}(0.4x_t+0.2x_{\eta t})z\widetilde N(dt,dz),\quad t\geq0,
	\end{aligned}
\end{equation}
where $x_0=1$, $\eta=0.5$ and $\nu\in\mathcal{P}(\mathbb{R}^d)$. The compensator is given by $\pi(dz)dt=2f(z)dudt$ where $f(z)\sim N(0,1)$. It can be also noticed that the coefficients of Eq.(\ref{example02}) satisfy Assumptions 2.1, 5.1 and meet the conditions of theorems 5.1 and 5.2 with $c_1=12, c_2, c_3=2$ and $c_4=1$ and $b_1=2, b_2=4, b_3=3$ and $b_4=6$. Also, $\varDelta t^{*}=\frac{(c_1-c_3)-(c_2+c_4)(1+[1/\eta])}{(b_1+b_3)+(b_2+b_4)(1+[1/\eta])}=0.114$. Therefore the solution of Eq.(\ref{example02}) possesses the property of almost sure exponential stability. Then $\lambda(r)=r^{5}$ and $\kappa(\varDelta t)=\varDelta t^{-1/4}$ are chosen so that the truncated Euler-Maruyama algorithm applied to the interacting particle system with respect to Eq.(\ref{example02}) is almost sure exponential stable. Moreover several trajectories of the numerical algorithm solutions with step size $\varDelta t=0.003$ where $\varDelta t\in(0,\varDelta t^{*})$ and $N=100$ particles to approximate the interacting particle system with respect to Eq.(\ref{example02}) are simulated and plotted in Figure \ref{fig:graph4} where it is clear that the average of these 100 particles' numerical solutions given by each trajectory is almost sure exponential stable.\\
\begin{figure}[!h]
	\centering
	\includegraphics[width=.5\textwidth]{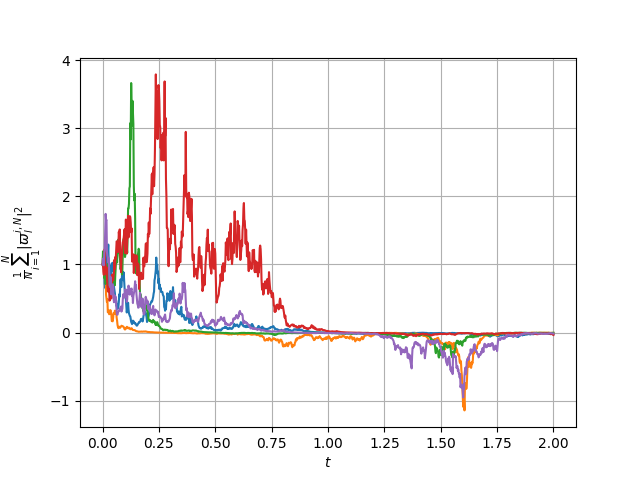}
	\caption{{\small  Trajectories of the numerical solutions of Eq.(\ref{example02}) over interval $[0,2]$.}}
	\label{fig:graph4}
\end{figure}\\}\\
{\small
\textbf{Funding:} There are no funders to report.}\\\\
{\small
\textbf{Conflicts of Interest:} This work does not have any conflicts of interest.}\\\\
{\small
\textbf{Data Availability Statement:} Data sharing not applicable to this article.}

\bibliographystyle{unsrt}
\bibliography{REFERENCES}

@article{gao2024convergence,
	title={Convergence rate in {L}p sense of tamed {E}{M} scheme for highly nonlinear neutral multiple-delay stochastic McKean--Vlasov equations},
	author={Gao, Shuaibin and Guo, Qian and Hu, Junhao and Yuan, Chenggui},
	journal={Journal of Computational and Applied Mathematics},
	volume={441},
	pages={115682},
	year={2024},
	publisher={Elsevier}
}

@article{mckean1966class,
	title={A class of {M}arkov processes associated with nonlinear parabolic equations},
	author={McKean Jr, Henry P},
	journal={Proceedings of the National Academy of Sciences},
	volume={56},
	number={6},
	pages={1907--1911},
	year={1966}
}

@article{biswas2020well,
	title={Well-posedness and tamed {E}uler schemes for {M}cKean-{V}lasov equations driven by L\'evy noise},
	author={Biswas, Sani and Kumar, Chaman and Reis, Gon{\c{c}}alo dos and Reisinger, Christoph and others},
	journal={arXiv preprint arXiv:2010.08585},
	year={2020}
}

@article{liu2023stability,
	title={Stability of the numerical scheme for stochastic {M}cKean-{V}lasov equations},
	author={Liu, Zhuoqi and Gao, Shuaibin and Yuan, Chenggui and Guo, Qian},
	journal={arXiv preprint arXiv:2312.12699},
	year={2023}
}

@article{liao1997almost,
	title={Almost sure exponential stability of neutral stochastic differential difference equations},
	author={Liao, XX and Mao, Xuerong},
	journal={Journal of Mathematical Analysis and Applications},
	volume={212},
	number={2},
	pages={554--570},
	year={1997},
	publisher={Elsevier}
}

@article{guhlke2018stochastic,
	title={Stochastic many-particle model for {LFP} electrodes},
	author={Guhlke, Clemens and Gajewski, Paul and Maurelli, Mario and Friz, Peter K and Dreyer, Wolfgang},
	journal={Continuum Mechanics and Thermodynamics},
	volume={30},
	pages={593--628},
	year={2018},
	publisher={Springer}
}

@article{dreyer2011phase,
	title={Phase transition in a rechargeable lithium battery},
	author={Dreyer, Wolfgang and Gaber{\v{s}}{\v{c}}ek, Miran and Guhlke, Clemens and Huth, Robert and Jamnik, Janko},
	journal={European Journal of Applied Mathematics},
	volume={22},
	number={3},
	pages={267--290},
	year={2011},
	publisher={Cambridge University Press}
}

@book{carmona2018probabilistic,
	title={Probabilistic theory of mean field games with applications I-II},
	author={Carmona, Ren{\'e} and Delarue, Fran{\c{c}}ois and others},
	year={2018},
	publisher={Springer}
}

@article{bossy2015clarification,
	title={Clarification and complement to “Mean-field description and propagation of chaos in networks of Hodgkin--Huxley and FitzHugh--Nagumo neurons”},
	author={Bossy, Mireille and Faugeras, Olivier and Talay, Denis},
	journal={The Journal of Mathematical Neuroscience (JMN)},
	volume={5},
	pages={1--23},
	year={2015},
	publisher={Springer}
}

@book{bensoussan2013mean,
	title={Mean field games and mean field type control theory},
	author={Bensoussan, Alain and Frehse, Jens and Yam, Phillip and others},
	volume={101},
	year={2013},
	publisher={Springer}
}

@article{wu2022stabilization,
	title={Stabilization of Stochastic {M}cKean-{V}lasov Equations with Feedback Control Based on Discrete-Time State Observation},
	author={Wu, Hao and Hu, Junhao and Gao, Shuaibin and Yuan, Chenggui},
	journal={SIAM Journal on Control and Optimization},
	volume={60},
	number={5},
	pages={2884--2901},
	year={2022},
	publisher={SIAM}
}

@article{benachour1998nonlinear,
	title={Nonlinear self-stabilizing processes--I Existence, invariant probability, propagation of chaos},
	author={Benachour, Sa{\i}d and Roynette, Bernard and Talay, Denis and Vallois, Pierre},
	journal={Stochastic processes and their applications},
	volume={75},
	number={2},
	pages={173--201},
	year={1998},
	publisher={Elsevier}
}

@article{bossy1997stochastic,
	title={A stochastic particle method for the {M}cKean-{V}lasov and the Burgers equation},
	author={Bossy, Mireille and Talay, Denis},
	journal={Mathematics of computation},
	volume={66},
	number={217},
	pages={157--192},
	year={1997}
}

@article{crisan2010approximate,
	title={Approximate {M}cKean-{V}lasov representations for a class of SPDEs},
	author={Crisan, Dan and Xiong, Jie},
	journal={Stochastics An International Journal of Probability and Stochastics Processes},
	volume={82},
	number={1},
	pages={53--68},
	year={2010},
	publisher={Taylor \& Francis}
}

@article{dos2022simulation,
	title={Simulation of {M}cKean-{V}lasov {S}{D}{E}s with super-linear growth},
	author={dos Reis, Gon{\c{c}}alo and Engelhardt, Stefan and Smith, Greig},
	journal={IMA Journal of Numerical Analysis},
	volume={42},
	number={1},
	pages={874--922},
	year={2022},
	publisher={Oxford University Press}
}

@article{guo2024convergence,
	title={Convergence analysis of an explicit method and its random batch approximation for the {M}cKean-{V}lasov equations with non-globally Lipschitz conditions},
	author={Guo, Qian and He, Jie and Li, Lei},
	journal={ESAIM: Mathematical Modelling and Numerical Analysis},
	volume={58},
	number={2},
	pages={639--671},
	year={2024},
	publisher={EDP Sciences}
}

@book{mao2006stochastic,
	title={Stochastic differential equations with {M}arkovian switching},
	author={Mao, Xuerong and Yuan, Chenggui},
	year={2006},
	publisher={Imperial college press}
}

@article{ockendon1971dynamics,
	title={The dynamics of a current collection system for an electric locomotive},
	author={Ockendon, John Richard and Tayler, Alan B},
	journal={Proceedings of the Royal Society of London. A. Mathematical and Physical Sciences},
	volume={322},
	number={1551},
	pages={447--468},
	year={1971},
	publisher={The Royal Society London}
}

@article{kou2002jump,
	title={A jump-diffusion model for option pricing},
	author={Kou, Steven G},
	journal={Management science},
	volume={48},
	number={8},
	pages={1086--1101},
	year={2002},
	publisher={INFORMS}
}

@article{svishchuk2000stochastic,
	title={The stochastic stability of interest rates with jump changes},
	author={Svishchuk, AV and Kalemanova, AV},
	journal={Theory of Probability and Mathematical Statistics},
	volume={48},
	number={61},
	pages={161--172},
	year={2000},
	publisher={Providence, American Mathematical Society, 1974-}
}

@book{tankov2003financial,
	title={Financial modelling with jump processes},
	author={Tankov, Peter},
	year={2003},
	publisher={Chapman and Hall/CRC}
}

@article{merton1976option,
	title={Option pricing when underlying stock returns are discontinuous},
	author={Merton, Robert C},
	journal={Journal of financial economics},
	volume={3},
	number={1-2},
	pages={125--144},
	year={1976},
	publisher={Elsevier}
}

@article{agrawal2020jump,
	title={Jump Models with Delay—Option Pricing and Logarithmic {E}uler--{M}aruyama Scheme},
	author={Agrawal, Nishant and Hu, Yaozhong},
	journal={Mathematics},
	volume={8},
	number={11},
	pages={1932--1953},
	year={2020},
	publisher={MDPI}
}

@article{maghsoodi1996mean,
	title={Mean square efficient numerical solution of jump-diffusion stochastic differential equations},
	author={Maghsoodi, Y},
	journal={Sankhy{\=a}: The Indian Journal of Statistics, Series A},
	pages={25--47},
	year={1996},
	publisher={JSTOR}
}

@article{higham2005numerical,
	title={Numerical methods for nonlinear stochastic differential equations with jumps},
	author={Higham, Desmond J and Kloeden, Peter E},
	journal={Numerische Mathematik},
	volume={101},
	number={1},
	pages={101--119},
	year={2005},
	publisher={Springer}
}

@article{bass2004stochastic,
	title={Stochastic differential equations driven by stable processes for which pathwise uniqueness fails},
	author={Bass, Richard F and Burdzy, Krzysztof and Chen, Zhen-Qing},
	journal={Stochastic processes and their applications},
	volume={111},
	number={1},
	pages={1--15},
	year={2004},
	publisher={Elsevier}
}

@article{hobson1998complete,
	title={Complete models with stochastic volatility},
	author={Hobson, David G and Rogers, Leonard CG},
	journal={Mathematical Finance},
	volume={8},
	number={1},
	pages={27--48},
	year={1998},
	publisher={Wiley Online Library}
}

@article{haghighi2019split,
	title={Split-step double balanced approximation methods for stiff stochastic differential equations},
	author={Haghighi, Amir and R{\"o}{\ss}ler, Andreas},
	journal={International Journal of Computer Mathematics},
	volume={96},
	number={5},
	pages={1030--1047},
	year={2019},
	publisher={Taylor \& Francis}
}

@book{milstein2004stochastic,
	title={Stochastic numerics for mathematical physics},
	author={Milstein, Grigori N and Tretyakov, Michael V},
	volume={39},
	year={2004},
	publisher={Springer}
}

@article{hutzenthaler2011strong,
	title={Strong and weak divergence in finite time of {E}uler's method for stochastic differential equations with non-globally {L}ipschitz continuous coefficients},
	author={Hutzenthaler, Martin and Jentzen, Arnulf and Kloeden, Peter E},
	journal={Proceedings of the Royal Society A: Mathematical, Physical and Engineering Sciences},
	volume={467},
	number={2130},
	pages={1563--1576},
	year={2011},
	publisher={The Royal Society Publishing}
}

@article{guo2018note,
	title={A note on the partially truncated Euler--Maruyama method},
	author={Guo, Qian and Liu, Wei and Mao, Xuerong},
	journal={Applied Numerical Mathematics},
	volume={130},
	pages={157--170},
	year={2018},
	publisher={Elsevier}
}

@article{hutzenthaler2012strong,
	title={Strong convergence of an explicit numerical method for {S}{D}{E}s with nonglobally {L}ipschitz continuous coefficients},
	author={Hutzenthaler, Martin and Jentzen, Arnulf and Kloeden, Peter E},
	journal={The Annals of Applied Probability},
	year={2012}
}

@article{tretyakov2013fundamental,
	title={A fundamental mean-square convergence theorem for {S}{D}{E}s with locally {L}ipschitz coefficients and its applications},
	author={Tretyakov, Michael V and Zhang, Zhongqiang},
	journal={SIAM Journal on Numerical Analysis},
	volume={51},
	number={6},
	pages={3135--3162},
	year={2013},
	publisher={SIAM}
}

@article{mao2015truncated,
	title={The truncated {E}uler--{M}aruyama method for stochastic differential equations},
	author={Mao, Xuerong},
	journal={Journal of Computational and Applied Mathematics},
	volume={290},
	pages={370--384},
	year={2015},
	publisher={Elsevier}
}

@article{mao2016convergence,
	title={Convergence rates of the truncated {E}uler--{M}aruyama method for stochastic differential equations},
	author={Mao, Xuerong},
	journal={Journal of Computational and Applied Mathematics},
	volume={296},
	pages={362--375},
	year={2016},
	publisher={Elsevier}
}

@article{guo2018truncated,
	title={The truncated {E}uler--{M}aruyama method for stochastic differential delay equations},
	author={Guo, Qian and Mao, Xuerong and Yue, Rongxian},
	journal={Numerical Algorithms},
	volume={78},
	pages={599--624},
	year={2018},
	publisher={Springer}
}

@article{mao2015existence,
	title={The existence and asymptotic estimations of solutions to stochastic pantograph equations with diffusion and {L}{\'e}vy jumps},
	author={Mao, Wei and Hu, Liangjian and Mao, Xuerong},
	journal={Applied Mathematics and Computation},
	volume={268},
	pages={883--896},
	year={2015},
	publisher={Elsevier}
}

@book{mao2007stochastic,
	title={Stochastic differential equations and applications},
	author={Mao, Xuerong},
	year={2007},
	publisher={Elsevier}
}

@article{bruti2007strong,
	title={Strong approximations of stochastic differential equations with jumps},
	author={Bruti-Liberati, Nicola and Platen, Eckhard},
	journal={Journal of Computational and Applied Mathematics},
	volume={205},
	number={2},
	pages={982--1001},
	year={2007},
	publisher={Elsevier}
}

@article{ahmadian2020exponential,
	title={Exponential mean-square stability of numerical solutions for stochastic delay integro-differential equations with {P}oisson jump},
	author={Ahmadian, Davood and Farkhondeh Rouz, Omid},
	journal={Journal of Inequalities and Applications},
	volume={2020},
	pages={1--33},
	year={2020},
	publisher={Springer}
}

@article{bruti2007approximation,
	title={Approximation of jump diffusions in finance and economics},
	author={Bruti-Liberati, Nicola and Platen, Eckhard},
	journal={Computational Economics},
	volume={29},
	pages={283--312},
	year={2007},
	publisher={Springer}
}

@article{albeverio2010existence,
	title={Existence of global solutions and invariant measures for stochastic differential equations driven by {P}oisson type noise with non-{L}ipschitz coefficients},
	author={Albeverio, Sergio and Brze{\'z}niak, Zdzis{\l}aw and Wu, Jiang-Lun},
	journal={Journal of Mathematical Analysis and Applications},
	volume={371},
	number={1},
	pages={309--322},
	year={2010},
	publisher={Elsevier}
}

@article{arriojas2007delayed,
	title={A delayed {B}lack and {S}choles formula},
	author={Arriojas, Mercedes and Hu, Yaozhong and Mohammed, Salah-Eldin and Pap, Gyula},
	journal={Stochastic Analysis and Applications},
	volume={25},
	number={2},
	pages={471--492},
	year={2007},
	publisher={Taylor \& Francis}
}

@article{meng2011pathwise,
	title={Pathwise estimation of stochastic differential equations with unbounded delay and its application to stochastic pantograph equations},
	author={Meng, Xuejing and Hu, Shigeng and Wu, Ping},
	journal={Acta Applicandae Mathematicae},
	volume={113},
	pages={231--246},
	year={2011},
	publisher={Springer}
}

@article{appleby2016sufficient,
	title={Sufficient conditions for polynomial asymptotic behaviour of the stochastic pantograph equation},
	author={Appleby, John AD and Buckwar, Evelyn},
	journal={arXiv preprint arXiv:1607.00423},
	year={2016}
}

@article{babasola2023stochastic,
	title={Stochastic delay differential equations: a comprehensive approach for understanding biosystems with application to disease modelling},
	author={Babasola, Oluwatosin and Omondi, Evans Otieno and Oshinubi, Kayode and Imbusi, Nancy Matendechere},
	journal={AppliedMath},
	volume={3},
	number={4},
	pages={702--721},
	year={2023},
	publisher={MDPI}
}

@article{vivek2022analysis,
	title={Analysis of stochastic pantograph differential equations with generalized derivative of arbitrary order},
	author={Vivek, Devaraj and Elsayed, Elsayed M and Kanagarajan, Kuppusamy},
	journal={e-Journal of Analysis and Applied Mathematics},
	volume={2022},
	number={1},
	year={2023},
	pages={24--32}
}

@article{ren2023stability,
	title={Stability and boundedness analysis of stochastic coupled systems with pantograph delay},
	author={Ren, Yong and Li, Jiaying},
	journal={International Journal of Control},
	volume={96},
	number={6},
	pages={1389--1396},
	year={2023},
	publisher={Taylor \& Francis}
}
\end{document}